\documentclass{imsart}

\usepackage{amsmath, amssymb, color}
\usepackage{amsthm} 
\usepackage[colorlinks,citecolor=blue,urlcolor=blue]{hyperref}
\usepackage{xcolor,colortbl}
\usepackage{mathabx}
\RequirePackage{enumerate}
\usepackage{color}
\usepackage{amsfonts, fancybox}
\usepackage{amssymb}
\usepackage[american]{babel}
\usepackage{psfrag,color}
\usepackage{dcolumn}
\usepackage[format=hang, justification=justified, singlelinecheck=false]{subcaption}
\usepackage{flafter, bm, dsfont}
\usepackage[section]{placeins}

\usepackage{multirow}
\usepackage{color}
\usepackage{amsmath}
\usepackage{float}
\usepackage{mathrsfs}
\usepackage{amsfonts}
\usepackage{dsfont}
\usepackage{amssymb}
\usepackage{algorithmic, algorithm}
\usepackage{wrapfig}

\usepackage{amssymb, latexsym}

\usepackage{graphicx}

\usepackage[square,sort,comma,numbers]{natbib}

\newtheorem{definition}{Definition}
\newtheorem{theorem}{Theorem}

\newtheorem{lemma}{Lemma}
\newtheorem{remark}{Remark}

\newtheorem{corollary}{Corollary}
\newtheorem{Assumption}{Assumption}
\newtheorem{proposition}{Proposition}
\newtheorem{openquestion}{Open question}

\newcommand{\R}{\mathbb R}
\newcommand{\PP}{\mathbb P}
\newcommand{\E}{\mathbb E}

\newcommand{\Q}{\mathbb Q}

\newcommand*\diff{\mathop{}\!\mathrm{d}}

\newcommand{\DS}{\displaystyle}

\begin{document}

\begin{frontmatter}

\title{Concentration of the empirical level sets of Tukey's halfspace depth}
\runtitle{Tukey depth level sets}

\begin{aug}

\author{\fnms{Victor-Emmanuel}~\snm{Brunel}\ead[label=veb]{vebrunel@math.mit.edu}},

\affiliation{Massachusetts Institute of Technology}

\address{{Victor-Emmanuel Brunel}\\
{Department of Mathematics} \\
{Massachusetts Institute of Technology}\\
{77 Massachusetts Avenue,}\\
{Cambridge, MA 02139-4307, USA}\\
\printead{veb}
}

\runauthor{Brunel, V.-E.}
\end{aug}

\begin{abstract}
Tukey's halfaspace depth has attracted much interest in data analysis, because it is a natural way of measuring the notion of depth relative to a cloud of points or, more generally, to a probability measure. Given an i.i.d. sample, we investigate the concentration of upper level sets of the Tukey depth relative to that sample around their population version. We show that under some mild assumptions on the underlying probability measure, concentration occurs at a parametric rate and we deduce moment inequalities at that same rate. In a computational prospective, we study the concentration of a discretized version of the empirical upper level sets.
\end{abstract}



\end{frontmatter}

\section{Preliminaries and notation}

\subsection{Preliminary}

Tukey's halfspace depth or, in short, Tukey depth, introduced by Tukey \citep{Tukey1975}, has attracted much attention in multivariate data analysis, as a tool for understanding and describing which data are relevant in a given cloud of points. For a finite multivariate sample, Tukey depth at any given point $x$ is the minimum proportion of points of the sample enclosed in a closed halfspace containing $x$. Tukey depth, together with other notions of statistical depths (see \citep{ZuoSerfling2000bis} for general definitions) has been studied and used extensively especially for description or graphical representation of data \citep{LiuPareliusSingh1999}, robust \citep{DonohoGasko1992, ArconesChenGine1994} or nonparametric (e.g., \citep{LiuSingh1993}) inference, bootstrap \citep{YehSingh1997}, supervised classification \citep{GhoshChaudhuri2005,GhoshChaudhuri2005bis}, etc... When the sample consists of i.i.d. random points, we call it empirical Tukey depth and it has a population analog (one can find formal definitions of the population Tukey depth in Euclidean spaces in \citep{RousseeuwRuts1999} and extensions to infinite dimensional Banach spaces in \citep{DuttaGhoshChaudhuri2011}). Consistency and limit theorems for the empirical Tukey depth are well-known (see \citep{Masse2004}, for instance, where the author tackles the asymptotic properties of the empirical Tukey depth seen as a stochastic process).

In this work, we are interested in the upper level sets of Tukey depth (we drop the qualifying \textit{upper} in the sequel). These sets are nested and the center of gravity of the deepest one is called the \textit{Tukey median}. On the opposite, the convex hull of a sample of $n$ points is the largest bounded empirical level set. Convergence and concentration of this random polytope has attracted a lot of attention in convex and stochastic geometry (see \citep{Brunel13''}, \citep{Fresen2013} and the references therein).

We show concentration of the level sets of the empirical Tukey depth of a given and fixed (independent of the sample size) level around the corresponding level sets of the population Tukey depth and we prove that the speed of convergence is parametric, i.e., of the order $n^{-1/2}$. Similar questions have already been tackled in earlier works. Consistency of the empirical level sets was proven in \citep{HeWang1997,ZuoSerfling2000} for general notions of statistical depth, including Tukey depth. In \citep{Kim2000}, the author shows that for all $\varepsilon\in (0,1)$, with probability $1-\varepsilon$, the empirical depth level set is sandwiched between two population level sets whose levels are at a distance of order $n^{-1/2}$ from each other. However, the constants are not explicit and the way they depend on $\varepsilon$ cannot be derived from the results, which, in turn, do not yield moment inequalities. In \citep{HeWang1997,ZuoSerfling2000,Kim2000}, the proofs are based on the global behavior of the stochastic process defined by the empirical depth, indexed by the ambient Euclidean space. Hence, the results in these works are based on global and very strong assumptions on the underlying probability measure. In our work, we focus on Tukey depth and only make local assumptions that guarantee some local continuity properties of the underlying distribution. We show that these assumptions are very weak, in the sense that they are satisfied by a broad class of distributions, including most commonly used ones. Not only we achieve the same (parametric) rate as obtained in \citep{Kim2000}, but our main result allows us to derive moment inequalities with a parametric rate.

Our approach is based on a polyhedral representation of the level sets of the population and empirical Tukey depths. As we will see in Lemma \ref{ThmEqualitySets}, which is a refinement of Theorem 2 in \citep{KongMizera2012}, these level sets can also be written as multivariate quantile sets, defined as convex regions that satisfy infinitely many linear constraints. It is because of such a multivariate quantile representation that the level sets of Tukey depth have also attracted attention in multivariate quantile regression (see \citep{Chaudhuri1996, HallinPaindaveineSiman2010} and the references therein). With this approach, we reduce the problem to that of estimating the support function of the population level sets. We believe that the techniques we use in our proofs could be useful in other problems related to support function estimation. For instance, in \citep{Guntuboyina2012}, the support function of an unknown convex set is observed up to some noise; We believe that our proof method could be used in order to bound from above the risk for estimation of the unknown convex set in Hausdorff distance, whereas the measure of the risk used in \citep{Guntuboyina2012} does not have a natural, geometric interpretation.

Computation of the empirical Tukey depth level sets for samples of $n$ points is a challenging problem. In dimension 2, they can be computed in $O(n^2)$ (see \citep{MillerRamaswamiRousseeuw2003}). A naive computation of the Tukey depth at one point would require to explore infinitely many halfspaces, which is not feasible. In higher dimensions, there is no practical and efficient way to compute the level sets of the Tukey depth. This is why we define a proxy for the empirical level sets, based on a discretized version of the Tukey depth. We show that they are consistent and still concentrate at the same parametric speed as the original ones. In practice, the number of operations required to compute this proxy grows exponentially with the dimension of the ambient space, but it can still be useful if the dimension is not too large.

Before going further into details, we introduce some notation. In this paper, $d\geq 2$ and $n\geq 1$ are fixed integers, unless stated otherwise. The Euclidean norm in $R^d$ is denoted by $|\cdot|$ and the dot product between two vectors $x$ and $y$ is denoted by $\langle x,y\rangle$. The $(d-1)$-dimensional unit sphere is $\mathcal S^{d-1}=\{u\in\R^d:|u|=1\}$. For $u\in\mathcal S^{d-1}$, $u^\perp$ stands for the hyperplane in $\R^d$ that is orthogonal to $u$. If $k$ is a positive integer, $a\in\R^k$ and $R\geq 0$, $B_k(a,R)$ (resp. $B_k'(a,R)$) stands for the closed (resp. open) Euclidean ball in $\R^k$ with center $a$ and radius $R$. When $k=d$, we drop the subscript $k$.

The complement of a set $A$ is denoted by $A^{\complement}$. The symmetric difference between two sets $A$ and $B$ in $\R^d$ is denoted by $A\triangle B$. For $k\geq 1$, if $A$ is a measurable set in $\R^k$ (equipped with the Lebesgue measure), we denote by $\textsf{Vol}_k(A)$ its $k$-dimensional volume, i.e., its Lebesgue measure in $\R^k$.

For $A\subseteq \R^d$, the interior of $A$ is denoted by $\overset{\circ}{A}$: this is the largest open set included in $A$. The collection of closed halfspaces in $\R^d$ is denoted by $\mathcal H$. For $u\in\mathcal S^{d-1}$ and $t\in\R$, we define the closed halfspace $\DS H_{u,t}=\{x\in\R^d:\langle u,x\rangle\leq t\}$.

The Hausdorff distance between two sets $K,K'\subseteq \R^d$ is 
	$$d_{\textsf{H}}(K,K')=\inf\{\varepsilon>0 : K\subseteq K'+\varepsilon B(0,1) \mbox{ and } K'\subseteq K+\varepsilon B(0,1)\},$$
where we set $\inf(\emptyset)=\infty$.
If $K$ is a convex body (i.e., convex and compact), its support function $h_K$ is defined as $\DS h_K(u)=\max_{x\in K} \langle u,x\rangle$, $u\in\R^d$.


The cardinality of a finite set $I$ is denoted by $\#I$. For $x\in\R$, we denote by $\lceil x \rceil$ the smallest integer larger or equal to $x$.

Throughout the paper, $X, X_1, X_2, \ldots$ are independent, identically distributed (i.i.d.) random variables defined on a probability space $(\Omega,\mathcal F,\PP)$, taking values in $\R^d$. Their common probability distribution is denoted by $\mu$ and is defined on the Borel $\sigma$-algebra of $\R^d$. The empirical distribution $\mu_n$ is defined by $\mu_n=\frac{1}{n}\sum_{i=1}^n \delta_{X_i}$, where $\delta_a$ is the Dirac measure at the point $a\in\R^d$.
 
For two positive sequences $(a_n)_{n\geq 1}$ and $(b_n)_{n\geq 1}$, we write $a_n=O(b_n)$ when the ratio $a_n/b_n$ is bounded uniformly in $n\geq 1$. For two positive sequences of random variables $(A_n)_{n\geq 1}$ and $(B_n)_{n\geq 1}$, we write $A_n=O_\PP(B_n)$ when for all $\delta>0$, there exists $M_\delta>0$ such that $\PP[A_n>M_\delta B_n]\leq\delta, \forall n\geq 1$.

Section \ref{SecGen} is devoted to general results about Tukey depth level sets. Our main theorems are given in Section \ref{SecMainResults} and the proofs are deferred to Section \ref{SecAppendix}. The rest of this section is dedicated to important definitions.

\subsection{Definitions}

The Tukey depth associated with a probability measure $\nu$ in $\R^d$ is the function $$D_{\nu}(x)=\inf_{H\in\mathcal H:x\in H}\nu(H), \hspace{3mm} \forall x\in\R^d.$$ 

We refer to $D_\mu$ as the \textit{population Tukey depth} and to $D_{\mu_n}$ as the \textit{empirical Tukey depth}.
 
In this work, we are interested in comparing the level sets of $D_\mu$ and $D_{\mu_n}$. Let $\alpha\in (0,1)$ be fixed. The $\alpha$-level set of $D_\mu$ is defined as $\DS G_{\mu}=\{x\in\R^d:D_{\mu}(x)\geq\alpha\}$ and we denote by $\hat G$ the $\alpha$-level set of $D_{\mu_n}$: $\DS \hat G=\left\{x\in\R^d:D_{\mu_n}(x)\geq \alpha\right\}$.
We study how fast $\hat G$ concentrates around $G_{\mu}$, i.e., how fast the stochastic convergence of $d_{\textsf{H}}(\hat G,G_{\mu})$ to zero is. 
As intermediate tools and for independent interest, we introduce the following sets associated with $\mu$:

\begin{enumerate}

	\item \textit{The multidimensional $(1-\alpha)$-quantile set of $\mu$:} \newline
	Let $X$ be a random variable with probability distribution $\mu$. For $u\in \R^d$, let $q_u^{\flat}$ and $q_u^{\sharp}$ be the lower and upper $(1-\alpha)$-quantile of $\langle u,X\rangle$, respectively:
	$$q_u^{\flat}=\inf\{t\in\R:\PP[\langle u,X\rangle\leq t]\geq 1-\alpha\} \quad \mbox{and}$$
	$$q_u^{\sharp}=\sup\{t\in\R:\PP[\langle u,X\rangle\geq t]\geq \alpha\}.$$
	The corresponding lower and upper multidimensional $(1-\alpha)$-quantile sets of $\mu$ are defined as
\begin{equation} \label{MultdimQuant}
G^{\eta}_{\textsf{MQ}}=\{x\in\R^d:\langle u,x\rangle\leq q_u^{\eta}, \forall u\in \mathcal S^{d-1}\}, \hspace{3mm} \eta\in\{\flat,\sharp\}.
\end{equation}
	
	\item \textit{The $\alpha$-floating body of $\mu$:} $\DS G_{\textsf{FB}}=\bigcap_{H\in\mathcal H:\mu(H)\geq 1-\alpha}H$.
\end{enumerate}

As we will see in Lemma \ref{ThmEqualitySets} below, these sets are other representations of the Tukey depth level sets. The representation in terms of multidimensional quantile sets is particularly convenient for our purposes because it characterizes the Tukey depth level sets through linear constraints.
We make the floating body part of our analysis because it plays an important role for random polytopes. Barany and Larman \citep{BaranyLarman1988} proved that if $\mu$ is the uniform distribution in a convex and compact set of volume 1, then the expected missing volume of the convex hull of $X_1,\ldots,X_n$ behaves aymptotically as the missing volume of the $(1/n)$-floating body of $\mu$. Fresen \citep{Fresen2013} proved that if $\mu$ is log-concave, the convex hull of $X_1,\ldots,X_n$ approximates the $(1/n)$-floating body of $\mu$ with high probability. For very small values of $\alpha$, even smaller than $1/n$, when the empirical level set would be a very poor estimator of $G_\mu$, \citep{He2016} defines and studies an estimator that extends univariate estimators from extreme value theory.

\section{Tukey depth level sets} \label{SecGen}

We start with a simple lemma that shows the relationships between the sets defined above: The Tukey depth level sets, the lower and upper multidimensional quantile sets and the floating bodies. This lemma is a refinement of Theorem 2 in \citep{KongMizera2012} but we include its proof at the end for the sake of completeness.

\begin{lemma} \label{ThmEqualitySets} 
$\DS G_{\textsf{FB}}=G^{\flat}_{\textsf{MQ}}\subseteq G^{\sharp}_{\textsf{MQ}}=G_\mu$.

\end{lemma}

In particular, if $\mu$ satisfies some continuity property, e.g., Assumption \ref{AssA} below, then $q_u^{\flat}=q_u^{\sharp}$ for all unit vectors $u$, so the inclusion becomes an equality and all four sets are equal. 

\citep{KongMizera2012} provides an interesting discussion about the multivariate quantile representation of $G_\mu$: In brief, the knowledge of $G_\mu$ does not imply the knowledge of all univariate quantiles $q_u^{\sharp}, u\in\mathcal S^{d-1}$. Indeed, some of the linear constraints that define $G^{\sharp}_{\textsf{MQ}}$ may not be active, i.e., there may be some unit vectors $u$ for which $\langle u,x\rangle<q_u^{\sharp}, \forall x\in G^{\sharp}_{\textsf{MQ}}$. This fact constitutes the main difficulty in the proof of Theorem \ref{MainTheorem} below, where we use the support function of $G^{\sharp}_{\textsf{MQ}}$. For $u\in\mathcal S^{d-1}$, it is clear that the linear constraint ``$\langle u,x\rangle\leq q_u^{\sharp}$" is active if and only if $h_{G_\mu}(u)=q_u^{\sharp}$. If that constraint is not active, then $h_{G_\mu}(u)<q_u^{\sharp}$. In that case, not only $G_\mu$ provides no information about $q_u^{\sharp}$, as discussed in \citep{KongMizera2012}, but $q_u^{\sharp}$ alone does not give any information about $h_{G_\mu}(u)$, and we need to understand how $h_{G_\mu}(u)$ depends on the $q_v^{\sharp}$'s that correspond to active constraints. 

For its independent interest, we may ask the following question: For which distributions $\mu$ are all the linear constraints that determine $G^{\sharp}_{\textsf{MQ}}$ active ? First, we have the following proposition about polyhedral representations of convex sets.

\begin{proposition} \label{LemmaPolyhedralRep}
	Let $(t_u)_{u\in\R^d}\subseteq \R$ be positively homogeneous, i.e., $t_{\lambda u}=\lambda t_u, \forall \lambda\geq 0, u\in\R^d$ and define the convex set $\DS G=\{x\in\R^d:\langle u,x\rangle \leq t_u, \forall u\in\mathcal S^{d-1}\}$. Assume that $\overset{\circ}{G} \neq \emptyset$. Then, the following statements are equivalent:
\begin{enumerate}[(i)]
	\item All the linear constraints that define $G$ are active;
	\item For all $u\in\mathcal S^{d-1}$, $h_G(u)=t_u$;
	\item The family $(t_u)_{u\in\R^d}$ is subadditive, i.e., $\DS t_{u+v}\leq t_u+t_v, \forall u,v\in\R^d$.
\end{enumerate}
\end{proposition}

As a consequence of this lemma, the upper quantiles $q_u^{\sharp}, u\in\mathcal S^{d-1}$, are completely determined by $G_\mu$ if and only if the family $\displaystyle{(q_u^{\sharp})_{u\in\R^d}}$ is sublinear, i.e., subadditive and positively homogeneous. 

\begin{openquestion} \label{OQ1}
	For what distributions $\mu$ are the upper quantiles $q_u^{\sharp}, u\in\R^d$, sublinear, no matter the value of $\alpha\in(0,1)$ ?
\end{openquestion}

A Gaussian distribution has sublinear upper quantiles, as a consequence of the triangle inequality for symmetric positive semidefinite matrices. If $\mu$ is the Gaussian distribution with centroid $m$ and covariance matrix $\Sigma$, then for all $u\in\mathcal S^{d-1}$, $q_u^{\sharp}=\langle u,m\rangle+\Phi^{-1}(1-\alpha)\sqrt{\Sigma(u,u)}$, where $\Phi$ is the cumulative distribution function of the univariate standard Gaussian distribution. The triangle inequality ensures that the map $u\in\R^d\mapsto \sqrt{\Sigma(u,u)}$ is sublinear, yielding sublinearity of $\displaystyle{(q_u^{\sharp})_{u\in\R^d}}$.

As a generalization of Gaussian distributions, and because they are known to be rigid (see \citep{LovaszVempala2007} for examples of this rigidity), we may ask if a log-concave probability measure have sublinear upper quantiles.

\begin{openquestion} \label{OQ2}
	Assume that $\mu$ is log-concave. Is it true that the upper quantiles $q_u^{\sharp}, u\in\R^d$ are sublinear, no matter the value of $\alpha\in(0,1)$ ?
\end{openquestion}

\begin{remark}
The multidimensional quantile sets are convex sets. Thus, they fail to capture the structure of complex probability measures, such as mixtures. The floating body (also called \textit{convex floating body} in the convex geometry literature, see \citep{SchuttWerner1990}) is defined as an intersection of closed halfspaces, i.e., the complement of the union of open halfspaces. Instead, one could think of an $r$-convex floating body, using the notion of $r$-convexity (see \citep{ManiLevitska1993}): $\DS G_{\textsf{FB}}^{(r)}=\left(\bigcup_{a\in\R^d:\mu(B'(a,r))<\alpha}B'(a,r)\right)^{{{\complement}}}$ and its empirical analog $\hat G_{\textsf{FB}}^{(r)}$ can be defined similarly, by replacing $\mu$ with $\mu_n$. When $r=\infty$, $G_{\textsf{FB}}^{(r)}=G_{\textsf{FB}}$. An asymptotic analysis of $\hat G_{\textsf{FB}}^{(r)}$ would require a different approach than ours, but seems to be relevant in order to describe more complex probability measures. In \citep{Pateiro2008}, $r$-convexity is exploited to estimate the support of probability distributions while relaxing convexity and even connectivity assumptions. We leave this question for further work.
\end{remark}

The next result shows that unless $\mu$ has atoms, the level set $G_\mu$ is empty when $\alpha$ is too large. 

\begin{theorem} \label{EmptyLevelSets}
	Let $\alpha>1/2$. Then, either $G_\mu$ is empty or it contains exactly one point. In the latter case, i.e., if $G_\mu=\{x\}$ for some $x\in\R^d$, then $x$ is an atom of $\mu$: $\mu(\{x\})>0$.
\end{theorem}

On the one hand, if $\mu$ has an atom $x$ with $\mu(\{x\})>1/2$, then $D_\mu(x)\geq \mu(\{x\}) >1/2$, hence, $G_\mu\neq\emptyset$ for $\alpha=\mu(\{x\})>1/2$. On the other hand, it is known (\citep{DonohoGasko1992}, Lemma 6.3) that $G_\mu$ is always nonempty when $\alpha\geq 1/(d+1)$. The following two examples show that very general probability measures $\mu$ can satisfy $G_\mu\neq\emptyset$ for large values of $\alpha\leq 1/2$, independent of the dimension $d$, and yet have no atoms:
\begin{itemize}
	\item If $\mu$ is centrally symmetric, i.e., satisfies $\mu(x+A)=\mu(x-A)$ for all Borel set $A\subseteq \R^d$, where $x$ is the center of symmetry of $\mu$, then $D_\mu(x)\geq 1/2$, hence, $G_\mu$ is nonempty for all $\alpha\in [0,1/2]$.
	\item If $\mu$ is log-concave, then any closed halfspace $H$ containing the centroid of $\mu$ satisfies $\mu(H)\geq e^{-1}$ (see Lemma 5.12 in \citep{LovaszVempala2007}). Hence, the depth of the centroid of $\mu$ is at least $e^{-1}$, which implies that $G_\mu$ is non empty for $\alpha$ as large as $e^{-1}\approx .37$.
\end{itemize}

\section{Concentration of the empirical Tukey depth level sets} \label{SecMainResults}

Consider the following assumptions, where we let $\varepsilon, L, r,R$ be fixed positive numbers satisfying $\varepsilon<r\leq R$

\begin{Assumption} \label{AssA}

\begin{itemize}

	\item For all $u\in\mathcal S^{d-1}$, the cumulative distribution function $F_u$ of $\langle u,X\rangle$ is continuous on $[q^{\sharp}_u-\varepsilon,q^{\sharp}_u+\varepsilon]$.
	
	\item $\DS F_u(t')-F_u(t)\geq L(t'-t)$, for all $u\in\mathcal S^{d-1}$ and all $t,t'\in\R$ with $q^{\sharp}_u-\varepsilon\leq t\leq t'\leq q^{\sharp}_u+\varepsilon$.
\end{itemize}

\end{Assumption}

\begin{Assumption} \label{AssB}

There exists $a\in\R^d$ such that $\DS B(a,r)\subseteq G_\mu \subseteq B(a,R)$.

\end{Assumption} 

Assumption \ref{AssA} ensures that $q_u^\flat=q_u^\sharp$ for all $u\in\mathcal S^{d-1}$, hence, that $G^{\flat}_{\textsf{MQ}}=G^{\sharp}_{\textsf{MQ}}$ and that the cumulative distribution functions $F_u$ are not too flat around their quantiles $q_u^{\flat}=q_u^{\sharp}$.


By Lemma \ref{ThmEqualitySets}, $\hat G$ can also be written as the empirical upper multidimensional $(1-\alpha)$-quantile set associated with $X_1,\ldots,X_n$: 

\begin{equation} \label{DefEmpMQ}
\hat G=\{x\in\R^d: \langle u,x\rangle\leq\hat q_u^{\sharp}, \forall u\in\mathcal S^{d-1}\},
\end{equation}
where, for $u\in\R^d$, $\DS \hat q_u^{\sharp}=\sup\Big\{t\in\R:\#\{i=1,\ldots,n:\langle u,X_i\rangle\geq t\}\geq n\alpha\Big\}$ is the upper empirical $(1-\alpha)$-quantile of $\langle u,X_1\rangle, \ldots, \langle u,X_n\rangle$. For the sake of notation, we will write $\hat q_u$ instead of $\hat q^{\sharp}_u$ in the sequel.

As a consequence of Lemma \ref{ThmEqualitySets}, in order to show concentration of $\hat G$ around $G_\mu$, one can compare their polyhedral representations given by \eqref{DefEmpMQ} and $G^{\sharp}_{\textsf{MQ}}$, which are written in terms of linear constraints. This is essential in the proof of our next theorem, which uses semi-infinite linear programming as one of its main ingredients.

Next theorem asserts that if Assumptions \ref{AssA} and \ref{AssB} are both satisfied, then $\hat G$ concentrates around $G_{\mu}$ at a parametric speed. In particular, that speed depends on the dimension $d$ only through multiplicative constants.

\begin{theorem} \label{MainTheorem}
Let $\mu$ satisfy Assumptions \ref{AssA} and \ref{AssB}. Then, the random set $\hat G$ satisfies the following deviation inequality:
\begin{equation*}
	\PP\left[d_{\textsf H}(\hat G,G_\mu)>\frac{Cx}{\sqrt n}\right]\leq Ae^{-L^2x^2/2+10\sqrt{5(d+1)}x},
\end{equation*}
for all $x\geq 0$ with $\DS \frac{10\sqrt{5(d+1)}}{L}\leq x<\varepsilon\sqrt n$, where $\displaystyle{C=\frac{R}{r}\frac{1+\varepsilon/r}{1-\varepsilon/r}}$ and $A=e^{-250(d+1)}$.
\end{theorem}

Note that in Theorem \ref{MainTheorem}, if $n$ is not large enough, the domain for $x$ will be empty. The constants depend on $d$ and the parameters $\varepsilon, r, R, L$. These parameters are hard to compute in practice, for a given distribution $\mu$. However, we give simple asymptotic consequences of Theorem \ref{MainTheorem} below. 

First, a truncated version of $\hat G$ has its expected error converging to zero at the speed $n^{-1/2}$:

\begin{corollary}\label{CorollFubini}

Define the random set 
$$\hat G^*= \begin{cases} \hat G\cap B'(0,\log n) \hspace{2mm} \mbox{if } \hat G\neq\emptyset\\
\{0\} \hspace{2mm} \mbox{otherwise.} \end{cases}$$ 
Let $\mu$ satisfy Assumptions \ref{AssA} and \ref{AssB} and assume, in addition, that $|a|\leq \tau$ for some $\tau>0$. Then, for all $k>0$, $\DS \E\left[d_{\textsf{H}}(\hat G^*,G_\mu)^k\right]=O\left(n^{-k/2}\right)$. The multiplicative constants in these asymptotic comparisons depend on $d,r,R,\varepsilon,L,\tau$ and $k$ only.

\end{corollary}

\begin{remark}
\begin{itemize}
	\item In Corollary \ref{CorollFubini}, the upper bounds are uniform on the class of probability measures $\mu$ that satisfy both Assumptions \ref{AssA} and \ref{AssB} with $|a|\leq \tau$. Hence, Corollary \ref{CorollFubini} gives an upper bound for the rate of the minimax risk in estimation of $G_\mu$ on that class of probability measures, and this rate is parametric. Note that the assumption $|a|\leq\tau$ could be dropped in Corollary \ref{CorollFubini}, but then the multiplicative constants in the asymptotic comparisons would also depend on $a$ and we would loose uniformity of the upper bounds.
	\item The threshold $\log n$ in the definition of $\tilde G$ is arbitrary and could be replaced with any sequence that grows to infinity at most polynomially in $n$.
\end{itemize}
\end{remark}

Define the maximal depth $\alpha_\mu^*$ of $\mu$ as $\displaystyle{\max_{x\in\R^d}D_\mu(x)}$. Consider the two following assumptions:

\begin{Assumption} \label{Ass11}
The probability measure $\mu$ is absolutely continuous with respect to the Lebesgue measure, its density $f$ is continuous and positive everywhere and there exist $C>0$ and $\nu>d-1$ such that $\displaystyle{|f(x)|\leq C\left(1+|x|\right)^{-\nu}}, \forall x\in\R^d$.
\end{Assumption}

In the sequel, if $\mu$ has a density $f$ with respect to the Lebesgue measure, we call the support of $\mu$ the set of vectors $x\in\R^d$ for which $f(x)>0$.

\begin{Assumption} \label{Ass12}
The probability measure $\mu$ is absolutely continuous with respect to the Lebesgue measure, its support is bounded and convex and its density is uniformly continuous on its support. 
\end{Assumption}


Assumptions \ref{Ass11} and \ref{Ass12} are sufficient but not necessary for next corollary. However, they include a lot of useful distributions. For example, any log-concave distribution in $\R^d$ with positive density satisfies Assumption \ref{Ass11}: A log-concave density is continous on its support and decays exponentially fast when $|x|\to\infty$. If $\mu$ has a density of the form $f(x)=h(\langle x,\Sigma x\rangle)$, where $\Sigma$ is a $d\times d$ symmetric positive definite matrix and $h$ is a positive continuous function that satisfies $h(t)\leq C(1+|t|)^{-\nu}$ for all $t\in\R$, with $\nu>d-1$, then $\mu$ satisfies Assumption \ref{Ass11} as well. If $\mu$ is the uniform distribution on a compact, convex set in $\R^d$, then it satisfies Assumption \ref{Ass12}.

\begin{corollary} \label{CorollaryStoch}
	Let $\mu$ satisfy either Assumption \ref{Ass11} or Assumption \ref{Ass12}. Suppose that $\alpha\in(0,\alpha_\mu^*)$, independently of $n$. Then, $\DS d_{\textsf{H}}(\hat G,G_\mu)=O_\PP\left(n^{-1/2}\right)$.
\end{corollary}

\begin{remark}
	\begin{itemize}
		\item Corollary \ref{CorollaryStoch} shows that the rate of convergence of the empirical level sets is parametric. 
		\item Surprisingly, if $\mu$ is the uniform distribution on a compact, convex set $K$ in $\R^d$, the rate does not depend on the smoothness of the boundary of $K$. This is paradoxical, since it is known that if $\alpha=1/n$, $\hat G$ is the convex hull of $X_1,\ldots,X_n$, which converges to $K$ at a rate that depends on the smoothness of the boundary of $K$ (see \citep{BaranyLarman1988}). However, in \citep{BaranyLarman1988}:
	\begin{itemize}
		\item $\alpha=1/n$ depends on $n$. In our work, $\alpha$ does not depend on $n$ and hence, the floating body $G_{\textsf{FB}}=G_\mu$ is bounded away from the boundary of $K$, which attenuates the effect of its smoothness.
		\item Convergence is towards the support $K$ itself, not towards the floating body of $\mu$. When $\alpha=1/n$, it is not clear whether the convergence of the distance between the empirical and the population $(1/n)$-convex bodies depends on the smoothness of the boundary of $K$. By the triangle inequality, $d_{\textsf H}(\hat G,K)\leq d_{\textsf H}(\hat G,G_\mu)+d_{\textsf H}(G_\mu,K)$. The $(1/n)$-floating body $G_\mu$ converges to $K$ at a speed that depends on the smoothness of the boundary of $K$ \citep{BaranyLarman1988,SchuttWerner1990}, but to the best of our knowledge, it is not known whether the speed of convergence of $d_{\textsf H}(\hat G,G_\mu)$ depends on the smoothness of $K$ too.
	\end{itemize}
	\item \citep{Kim2000} obtained the parametric rate $n^{-1/2}$ for general measures of statistical depth, under quite strong assumptions on $\mu$ which rule out many important distributions, as compared to ours (e.g., compactly supported densities). In addition, they do not compare $\hat G$ to $G_\mu$ directly, but to level sets of $D_\mu$ with levels $\alpha\pm Mn^{-1/2}$, for some $M>0$, leaving out a deterministic bias. Yet, we believe that they could achieve the same rate as ours. However, unlike Theorem \ref{MainTheorem}, their result is not informative about the tail of the distribution of $d_{\textsf H}(\hat G,G_\mu)$, because of implicit dependency of the constant $M$ on the probability level (see \citep{Kim2000}, Theorem 1).
\end{itemize}
\end{remark}

Computation of $\hat G$ is a hard problem. Its concentration around $G_\mu$ is a question of its own geometric and probabilistic interest, but it also has important statistical implications. For instance, as we saw in Corollary \ref{CorollFubini}, it provides a benchmark for the minimax risk for estimation of $G_\mu$ based on an i.i.d. sample. However, if $\hat G$ is too hard to compute, this does not have much of a practical interest. 
Computation of the Tukey depth $D_{\mu_n}$ at a single point is equivalent to the problem of finding a hemisphere that contains the largest number of points positioned on the unit sphere, which is NP hard in high dimension \citep{JohnsonPreparata1978}. However, in fixed dimension, some deterministic and random algorithms to compute an approximate or exact value of the Tukey depth have been suggested (see \citep{RousseeuwRuts1996,RousseeuwStruyf1998,DyckerhoffMozharovskyi2016}  and the references therein).
For the actual computation of the Tukey depth level sets relative to a point cloud in dimension 2, we refer to \citep{MillerRamaswamiRousseeuw2003}. These sets are polygons, hence, their computation reduces to finding either their vertices or their faces. To our knowledge, there are no algorithms to compute these sets exactly when $d\geq 3$. Here, we define a random approximation of $\hat G$ that can be computed exactly, yet in an exponential time in $d$. Lemma \ref{ThmEqualitySets} gives a representation of $\hat G$ through infinitely many linear constraints. By selecting a finite number of these constraints, using a collection of unit vectors that are well spread on the unit sphere, one can obtain a suitable approximation of $\hat G$. 

Our random approximation is obtained by sampling random vectors on the unit sphere. If $M$ is a positive integer, denote by $\DS \tilde G_M=\left\{x\in\R^d: \langle U_j,x\rangle \leq \hat q_{U_j}, \hspace{2mm} \forall j=1,\ldots,M\right\}$, where $U_1,\ldots,U_M$ are i.i.d. uniform random variables on $\mathcal S^{d-1}$, independent of $X_1,\ldots,X_n$. The following theorem shows that a certain choice of $M$ leads to an estimator that of $G_\mu$ that concentrates as fast as $\hat G$.

\begin{theorem} \label{TheoremComput}

Let $\mu$ satisfy Assumptions \ref{AssA} and \ref{AssB} and assume that the quantiles $\DS (q_u^\sharp)_{u\in\R^d}$ are subadditive. Then, for all $M\geq 1$, the random set $\tilde G_M$ satisfies the following deviation inequality:
\begin{align*} 
	& \PP\left[d_{\textsf{H}}(\tilde G_M,G_\mu) > \frac{Cx+4R}{\sqrt n}\right] \\
	& \hspace{10mm} \leq Ae^{-L^2x^2/2+10\sqrt{5(d+1)}Lx}+6^d\exp\left(-\frac{M}{2d8^{(d-1)/2}n^{d-1}}+(d/2)\log n\right),
\end{align*}
for all real number $x$ with $\DS 10\sqrt{5(d+1)}\leq x<\varepsilon\sqrt n$, where $c'_1, c'_2$ and $c'_3$ are positive constants that depend on $d, r, R, \varepsilon, L$ and $\rho$ only.

\end{theorem}

The explicit values of the constants can be easily derived from the proof.

\begin{remark}
	In Theorem \ref{TheoremComput}, we assume that the population quantiles are subadditive, which, by Proposition \ref{LemmaPolyhedralRep}, ensures that they are completely characterized by the knowledge of $G_\mu$. How strong this assumption is is an open question (see Open questions \ref{OQ1} and \ref{OQ2}).
\end{remark}

Theorem \ref{TheoremComput} yields the following asymptotic upper bound for a truncated version of $\tilde G_M$, if $M$ is chosen large enough.

\begin{corollary}\label{CorollFubiniComput}

Define the random set $\tilde G_M^*$ as
$$\hat G_M^*= \begin{cases} \tilde G_M\cap B(0,\log n) \hspace{2mm} \mbox{if } \tilde G_M\neq\emptyset, \\
\{0\} \hspace{2mm} \mbox{otherwise.} \end{cases}$$
Let $k>0$. Recall the notation and assumptions of Theorem \ref{TheoremComput}. If, in addition, $|a|\leq \tau$ for some $\tau>0$, then for $\DS M>2d8^{(d-1)/2}\frac{d+k}{2}n^{d-1}\log n$, $\tilde G_M^*$ satisfies $\DS \E\left[d_{\textsf{H}}(\tilde G_M^*,G_\mu)^k\right]=O\left(n^{-k/2}\right)$. The multiplicative constants in this asymptotic comparison depend on $d,r,R,\varepsilon,L$ and $k$ only.

\end{corollary}

In addition, the following stochastic upper bound holds under subadditivity of the population quantiles and either Assumption \ref{Ass11} or Assumption \ref{Ass12} :

\begin{corollary}
	Let $\alpha\in(0,\alpha_\mu^*)$ and $\mu$ satisfy either Assumption \ref{Ass11} or Assumption \ref{Ass12}. Then, if the quantiles $\DS (q_u^\sharp)_{u\in\R^d}$ are subadditive, $\DS d_{\textsf{H}}(\tilde G,G_\mu)=O_\PP\left(1/\sqrt n\right)$.
\end{corollary}

\section{Proofs}\label{SecAppendix}

\subsection{Preliminary lemmas in convex geometry and semi-infinite linear programming}

\begin{lemma}\label{LemmaBounded}

Let $K,L$ be two convex sets. Then, $\DS K\subseteq L \iff h_K(u)\leq h_L(u), \forall u\in\mathcal S^{d-1}$.
In particular, $K$ is bounded if and only if the restriction of its support function to the unit sphere is bounded.

\end{lemma}

\paragraph{Proof: }

The first part of the lemma follows from directly from the definition of the support function. For the second part, note that the support function of a ball centered at the origin with radius $R\geq 0$ is constant, equal to $R$ on the unit sphere. Hence, $\DS h_K(u)\leq R, \forall u\in \mathcal S^{d-1} \iff K\subseteq B(0,R)$, which proves the second part of the lemma. \hfill \textsquare

\vspace{4mm}

In the next two lemmas, we let $\phi:\mathcal S^{d-1}\to\R$ and $\DS K=\{x\in\R^d:\langle u,x\rangle \leq t_u, \hspace{4mm} \forall u\in\mathcal S^{d-1}\}.$

\begin{lemma} \label{Lemmacb}
The set $K$ is convex and compact.
\end{lemma}

\paragraph{Proof:} If $K$ is empty, then it is convex and compact. Assume that $K$ is nonempty. It is closed and convex, as the intersection of closed halfspaces. Let us show that $K$ is bounded, which will end the proof. Since $h_K$ is convex, it is continuous on the interior of its domain (the domain of $h_K$ is $\{u\in\R^d:h_K(u)<\infty\}$. For all $u\in\mathcal S^{d-1}$, $h_K(u)\leq t_u<\infty$, yielding that $h_K$ is continuous on $\mathcal S^{d-1}$. Since $\mathcal S^{d-1}$ is compact, the restriction of $h_K$ on the sphere needs to be bounded. Hence, by Lemma \ref{LemmaBounded}, $K$ is bounded. \hfill \textsquare

\begin{lemma} \label{LemmaInterior}
	If $\phi$ is continuous and $x\in\R^d$, then $\DS x\in \overset{\circ}{K} \iff \langle u,x\rangle < \phi(u), \forall u\in\mathcal S^{d-1}$.
\end{lemma}

\paragraph{Proof:} Let $x\in \overset{\circ}{K}$. Then, $B'(x,\eta)\subseteq K$ for some $\eta>0$. Let $u\in\mathcal S^{d-1}$. Then, $x+\eta u\in K$, yielding $\langle u,x+\eta u\rangle\leq\phi(u)$. Hence, $\langle u,x\rangle\leq \phi(u)-\eta<\phi(u)$ and this has to be true for all $u\in\mathcal S^{d-1}$. Now, let $x\in\R^d$ satisfying $\langle u,x\rangle < \phi(u), \forall u\in\mathcal S^{d-1}$. The map $u\in\mathcal S^{d-1}\mapsto \phi(u)-\langle u,x\rangle$ is continuous and positive on the compact $\mathcal S^{d-1}$, hence, there exists $\eta>0$ such that for all $u\in\mathcal S^{d-1}$, $\phi(u)-\langle u,x\rangle\geq \eta$. Then, it is easy to verify that $B'(x,\eta)\subseteq K$, yielding $x\in \overset{\circ}{K}$. \hfill \textsquare

%
%

When a convex set is defined through a collection of linear inequalities indexed by the unit sphere, the support function at a given unit $u_0$ vector can be interpreted as the value of a semi-infinite linear program. The following lemma states that under a continuity assumption, $u_0$ needs to lie in the convex cone spanned by the constraints that are active at a point $x^*$ that is a solution of that linear program. Note that when the number of linear constraints is infine, the existence of active constraints is not granted, as the following example shows.

Let $u_0\in\mathcal S^{d-1}$ and $\DS G=\{x\in\R^d:\langle u,x\rangle \leq 1, \forall u\in\mathcal S^{d-1}\setminus\{u_0\},\langle u_0,x\rangle \leq 2\}$. Then, since it is also true that $G=B'(0,1)$, the value of the semi-infinite linear program $\DS \max\{\langle u_0,x\rangle: x\in G\}$ is 1, uniquely attained at $x^*=u_0$. Yet, no constraint is active at $x^*$.

\begin{lemma}\label{LemmaActiveCone}

	Let $\phi$ be a continuous function on $\mathcal S^{d-1}$ and let $\DS K=\{x\in\R^d:\langle u,x\rangle\leq\phi(u), \forall u\in\mathcal S^{d-1}\}$. Assume that $\overset{\circ}{K}\neq\emptyset$.
For all $u_0\in\mathcal S^{d-1}$, there exists $x^*\in K$ such that $h_K(u_0)=\langle u_0,x^*\rangle$. Moreover, there exists $I\subseteq \mathcal S^{d-1}$ such that
\begin{itemize}
	\item $\# I\leq d$,
	\item $\langle u,x^*\rangle=\phi(u), \forall u\in I$,
	\item $\DS u_0=\sum_{u\in I}\lambda_u u$, for some nonnegative numbers $\lambda_u, u\in I$.
\end{itemize}

\end{lemma}

\paragraph{Proof: }

By Lemma \ref{Lemmacb}, $K$ is compact, which grants the existence of $x^*$, since $K\neq\emptyset$. Let $I^*=\{u\in\mathcal S^{d-1}:\langle u,x^*\rangle = \phi(u)\}$ be the set of active constraints at $x^*$ and let us prove that $I^*$ is not empty. The rest will follow using Theorem 2 in \citep{LopezStill2007} (Slater's condition is satisfied since we assume that $K$ has nonempty interior).

If $I^*$ was empty, then 
\begin{equation}\label{eqnPos475}
\forall u\in\mathcal S^{d-1}, \langle u,x^*\rangle < \phi(u).
\end{equation} 
Since the function $u\in\mathcal S^{d-1}\mapsto \phi(u)-\langle u,x^*\rangle $ is continuous and positive on the compact $\mathcal S^{d-1}$, there is a positive number $\eta$ such that $\phi(u)-\langle u,x^*\rangle\geq \eta, \forall u\in\mathcal S^{d-1}$. Hence, for all $u\in\mathcal S^{d-1}$
\begin{align*}
	\langle u,x^*+\eta u_0\rangle & = \langle u,x^*\rangle +\eta \langle u,u_0\rangle \\
	& \leq \phi(u)-\eta+\eta \langle u,u_0\rangle \\
	& \leq \phi(u)-\eta+\eta =\phi(u),
\end{align*}
yielding that $x^*+\eta u_0\in K$. This contradicts the maximality of $h_K(u_0)$, since $\langle u_0,x^*+\eta u_0\rangle > \langle u_0,x^*\rangle = h_K(u_0)$. \hfill \textsquare

\vspace{3mm}
In the next two lemmas, for any map $\zeta:\mathcal S^{d-1}\to\R$, we denote by $\DS G_\zeta=\{x\in\R^d:\langle u,x\rangle \leq \zeta(u), \forall u\in\mathcal S^{d-1}\}$
\begin{lemma} \label{Lemma_q_to_h}
	Let $\phi$ and $\hat\phi$ be two continuous functions on $\mathcal S^{d-1}$. Assume that $G_\phi$ and $G_{\hat\phi}$ have nonempty interiors. Let $R>r>0$ and assume that $B'(0,r)\subseteq G_\phi\subseteq B'(0,R)$. Let $\eta=\max_{u\in\mathcal S^{d-1}}|\hat\phi(u)-\phi(u)|$. If $\eta<r$, then $\DS d_{\textsf H}(G_{\hat\phi},G_\phi)\leq \frac{\eta R}{r}\frac{1+\eta/r}{1-\eta/r}$.
\end{lemma}

\paragraph{Proof: }

Let $u_0\in\mathcal S^{d-1}$. By Lemma \ref{LemmaActiveCone}, there exist $x\in G_\phi, \hat x\in G_{\hat\phi}$, $I,\hat I\subseteq \mathcal S^{d-1}$ with $\# I\leq d, \# {\hat I}\leq d$, such that $h_{G_\phi}(u_0)=\langle u_0,x\rangle$, $h_{G_{\hat\phi}}(u_0)=\langle u_0,{\hat x}\rangle$, $\langle u,x\rangle=\phi(u), \forall u\in I$, $\langle v,{\hat x}\rangle={\hat\phi}(v), \forall v\in {\hat I}$ and $u_0=\sum_{u\in I}\lambda_u u=\sum_{v\in {\hat I}}\hat\lambda_v v$, for some nonnegative families $\DS \left(\lambda_u\right)_{u\in I},\left(\hat\lambda_v\right)_{v\in \hat I}\geq 0$. Note that necessarily, for all $u\in I$ and $v\in\hat I$, $\phi(u)=h_{G_\phi}(u)$ and ${\hat\phi}(v)=h_{G_{\hat\phi}}(v)$. Then,
\begin{align}
	h_{G_{\hat\phi}}(u_0) & = h_{G_{\hat\phi}}\left(\sum_{u\in I}\lambda_u u\right) \leq \sum_{u\in I}\lambda_u h_{G_{\hat\phi}}(u) \leq \sum_{u\in I}\lambda_u {\hat\phi}(u) \leq \sum_{u\in I}\lambda_u (\phi(u)+\eta) \nonumber \nonumber \\
	\label{qtoh1} & = \sum_{u\in I}\lambda_u\langle u,x\rangle +\eta \sum_{u\in I}\lambda_u = \langle u_0,x\rangle+\eta \sum_{u\in I}\lambda_u = h_{G_\phi}(u_0)+\eta \sum_{u\in I}\lambda_u. 
\end{align}
In a similar fashion, we have that 
\begin{equation} \label{qtoh2}
	h_{G_\phi}(u_0) \leq h_{G_{\hat\phi}}(u_0)+\eta \sum_{v\in {\hat I}}\hat\lambda_v.
\end{equation}

By Lemma \ref{LemmaBounded} and since $B'(0,r)\subseteq G_\phi\subseteq B'(0,R)$, $r\leq h_{G_\phi}(u)\leq R$, for all $u\in\mathcal S^{d-1}$, yielding $\DS R\geq \langle u_0,x\rangle = \sum_{u\in I}\lambda_u\langle u,x\rangle = \sum_{u\in I}\lambda_u h_{G_\phi}(u)\geq r\sum_{u\in I}\lambda_u$. Hence, $\DS \sum_{u\in I}\lambda_u\leq \frac{R}{r}$ and by \eqref{qtoh1},
\begin{equation} \label{qtoh3}
	h_{G_{\hat\phi}}(u_0)\leq h_{G_\phi}(u_0)+\frac{\eta R}{r}.
\end{equation}

On the other hand,
\begin{equation} \label{qtoh4}
	\sum_{v\in {\hat I}}\hat\lambda_v\langle u,{\hat x}\rangle = \langle u_0,{\hat x}\rangle = h_{G_{\hat\phi}}(u_0)\leq h_{G_\phi}(u_0)+\frac{\eta R}{r}\leq R+\frac{\eta R}{r},
\end{equation}
where the third inequality comes from \eqref{qtoh3}. In addition, $\DS \sum_{v\in {\hat I}}\hat\lambda_v\langle u,{\hat x}\rangle = \sum_{v\in {\hat I}}\hat\lambda_v{\hat\phi}(v) \geq \sum_{v\in {\hat I}}\hat\lambda_v(\phi(v)-\eta)\geq \sum_{v\in {\hat I}}\hat\lambda_v(r-\eta)$,yielding, together with \eqref{qtoh4},
\begin{equation} \label{qtoh5}
	\sum_{v\in {\hat I}}\hat\lambda_v\leq \frac{R}{r}\frac{1+\eta/r}{1-\eta/r}.
\end{equation}
Finally, \eqref{qtoh1}, \eqref{qtoh2} and \eqref{qtoh5} yield
\begin{equation} \label{qtoh6}
	|h_{G_\phi}(u_0)-h_{G_{\hat\phi}}(u_0)| \leq \frac{R}{r}\frac{1+\eta/r}{1-\eta/r}.
\end{equation}
Since \eqref{qtoh6} is true for any arbitrary $u_0\in\mathcal S^{d-1}$, Lemma \ref{Lemma_q_to_h} is proven. \hfill \textsquare

\begin{definition} \label{DefNets}
	Let $\delta>0$. A $\delta$-net of the sphere $S^{d-1}$ is a subset $\mathcal N\subseteq \mathcal S^{d-1}$ such that $\DS \sup_{u\in\mathcal S^{d-1}}\inf_{v\in\mathcal N}|u-v|\leq\delta$.

\end{definition}

\begin{lemma} \label{Lemma_q_to_h_Disc}
Let $\delta\in (0,1)$ and $\mathcal N$ be a $\delta$-net of $\mathcal S^{d-1}$. Let $\phi$ and $\hat\phi:\R^d\to\R$, and assume that $\phi$ is sublinear. Let $r<R$ be two positive numbers and assume that $\DS B'(0,r)\subseteq G_\phi\subseteq B'(0,R)$. Let $\eta=\max_{u\in\mathcal N}|\phi(u)-\hat\phi(u)|$. If $\eta<r$, then $\DS d_{\textsf H}(G_\phi,G_{\hat\phi}^{\mathcal N})\leq \frac{\eta R}{r}\frac{1+\eta/r}{1-\eta/r}+\frac{2R\delta}{1-\delta}$.
\end{lemma}

\paragraph{Proof: } 

Before starting the proof, let us recall the following important property for support functions. If $K\subseteq B'(0,M)$ is a convex set, with $M>0$, then its support function is $M$-Lipschitz.

Let $G_{\phi}^{\mathcal N}=\{x\in\R^d:\langle u,x\rangle \leq \phi(u), \forall u\in\mathcal N\}$. By the triangle inequality,
\begin{equation} \label{qtohDisc1}
	d_{\textsf H}(G_\phi,G_{\hat\phi}^{\mathcal N})\leq d_{\textsf H}(G_\phi,G_{\phi}^{\mathcal N})+d_{\textsf H}(G_{\phi}^{\mathcal N},G_{\hat\phi}^{\mathcal N}).
\end{equation}

By Proposition \ref{LemmaPolyhedralRep}, $\phi(u)=h_{G_\phi}(u), \forall u\in\mathcal S^{d-1}$. Hence, $\phi(u)\leq R, \forall u\in\mathcal S^{d-1}$. Let $x\in G_{\phi}^{\mathcal N}$ with $x\neq 0$ and let $u=x/|x|$. Then, $|u-u^*|\leq\delta$ for some $u^*\in\mathcal N$, yielding $\DS |x|=\langle u,x\rangle = \langle u^*,x\rangle+\langle u-u^*,x\rangle \leq \phi(u^*)+\delta|x|\leq R+\delta |x|$. Hence, $\DS |x|\leq \frac{R}{1-\delta}$ and $G_{\phi}^{\mathcal N}\subseteq B'(0,R/(1-\delta))$. This entails that $h_{G_{\phi}^{\mathcal N}}$ is $R/(1-\delta)$-Lipschitz. Now, let $u_0\in\mathcal S^{d-1}$. On the one hand, since $G_{\phi}\subseteq G_{\phi}^{\mathcal N}$, $h_{G_{\phi}}(u_0)\leq h_{G_{\phi}^{\mathcal N}}(u_0)$. On the other hand, if $u^*\in\mathcal N$ satisfies $|u_0-u^*|\leq\delta$, then
\begin{align*}
	h_{G_{\phi}^{\mathcal N}}(u_0) & \leq h_{G_{\phi}^{\mathcal N}}(u^*)+\frac{R\delta}{1-\delta}\leq \phi(u^*)+\frac{R\delta}{1-\delta}=h_{G_{\phi}}(u^*)+\frac{R\delta}{1-\delta} \\
	& \leq h_{G_{\phi}}(u_0)+R|u_0-u^*|+\frac{R\delta}{1-\delta}\leq h_{G_{\phi}}(u_0)+R\delta+\frac{R\delta}{1-\delta} \\ 
	& = h_{G_{\phi}}(u_0)+\frac{2R\delta}{1-\delta} \leq h_{G_{\phi}}(u_0)+\frac{R\delta(2-\delta)}{1-\delta},
\end{align*}
where we used the fact that $h_{G_{\phi}}$ is $R$-Lipschitz. Therefore,
\begin{equation} \label{qtohDisc2}
	d_{\textsf H}(G_\phi,G_{\hat\phi}^{\mathcal N}) \leq \frac{2R\delta}{1-\delta}.
\end{equation}

Since $B'(0,r)\subseteq G_\phi\subseteq G_{\phi}^{\mathcal N}$, $G_{\phi}^{\mathcal N}$ has nonempty interior. So does $G_{\hat\phi}^{\mathcal N}$, since it is clear that $B'(0,r-\eta)\subseteq G_{\hat\phi}^{\mathcal N}$, using the facts that $\phi(u)\geq r, \forall u\in\mathcal S^{d-1}$, by Lemma \ref{LemmaBounded} and that $\eta<r$. Hence, using similar arguments as in the proof of Lemma \ref{Lemma_q_to_h},
\begin{equation} \label{qtohDisc3}
	d_{\textsf H}(G_{\phi}^{\mathcal N},G_{\hat\phi}^{\mathcal N})\leq \frac{\eta R}{r}\frac{1+\eta/r}{1-\eta/r}.
\end{equation}
Thus, \eqref{qtohDisc1}, \eqref{qtohDisc2} and \eqref{qtohDisc3} yield the desired result. \hfill \textsquare

\begin{lemma} \label{LemmaInteriorA}
	Let $K\subseteq \R^d$ be a convex set with nonempty interior. Let $A=\{(u,t)\in\mathcal S^{d-1}\times\R:(tu+u^\perp)\cap K\neq\emptyset\}$. Then, a pair $(u,t)\in\mathcal S^{d-1}\times\R$ is in $\overset{\circ}{A}$ if and only if there exists $\eta>0$ satisfying
\begin{equation} \label{EqLemmaInteriorA}
	(su+u^\perp)\cap \overset{\circ}{K}\neq\emptyset, \quad \forall s\in [t-\eta,t+\eta].
\end{equation}
\end{lemma}

\paragraph{Proof: }

Let $(u,t)\in\mathcal S^{d-1}\times\R$. 

Assume that $(u,t)\in\overset{\circ}{A}$. Then, there exists $\eta>0$ such that $(u,s)\in A$, for all $s\in [t-2\eta,t+2\eta]$. Let $s\in [t-2\eta,t+2\eta]$. Since $(u,s)\in A$, the affine hyperplane $su+u^\perp$ intersects $K$. It actually needs to intersect $\overset{\circ}{K}$. Indeed, $\overset{\circ}{K}$ is also the relative interior of $K$, since $K$ has nonempty interior. Hence, for the affine hyperplane $su+u^\perp$ to intersect $K$ but not its interior, it has to be a supporting hyperplane of $K$. This contradicts the fact that $K$ has elements on both sides of $su+u^\perp$.

Now, assume that $(u,t)\in\overset{\circ}{A}$ for some $\eta>0$. Then, the affine hyperplane $tu+u^\perp$ intersects the interior of $K$ and let $x\in (tu+u^\perp)\cap \overset{\circ}{K}$. Let $\eta>0$ be such that $B'(x,\eta)\subseteq \overset{\circ}{K}$. Let $\delta=\eta/(1+|x|)$ and $(v,s)\in\mathcal S^{d-1}\times\R$ with both $|v-u|\leq\delta$ and $|s-t|\leq\delta$. Since $x\in \langle v,x\rangle v+v^\perp$, the affine hyperplane $sv+v^\perp$ intersects $B'(0,x)$ if and only if $|s-\langle v,x\rangle|\leq\eta$, which holds by our choice of $\delta$. \hfill \textsquare

\begin{lemma} \label{Covariogram}
	Let $k$ be a positive integer and $K$ be a compact and convex set in $\R^k$ such that $0\in\overset{\circ}{K}$. Let $u\in\mathcal S^{k-1}$ and let $(u_n)_{n\geq 1}$ a sequence of unit vectors in $\R^k$ that converges to $u$. Let $(x_n)_{n\geq 1}$ be a sequence in $\R^k$ that converges to zero and $(U_n)_{n\geq 1}$ be a sequence of isometries in $\R^k$ that converges to the identity. Then, as $n\to 0$,

\begin{enumerate}
	\item $\displaystyle{\textsf{Vol}_{k-1}\left(\left((K+x_n)\cap u_n^\perp\right)\triangle (K\cap u_n^\perp)\right) \longrightarrow 0}$;
	\item $\displaystyle{\textsf{Vol}_{k-1}\left(\left(U_n(K)\cap u_n^\perp\right)\triangle (K\cap u_n^\perp)\right) \longrightarrow 0}$.
\end{enumerate}

\end{lemma}

\paragraph{Proof:} 

For $u\in\mathcal S^{k-1}$, set $\DS p_K(u)=\max\{\lambda\geq 0 : \lambda u\in K\}$. This is the (multiplicative) inverse of the gauge of $K$. Since $0\in\overset{\circ}{K}$, there exists $m>0$ such that $B_k'(0,m)\subseteq K$, yielding $p_K(u)\geq m$ for all $u\in\mathcal S^{k-1}$.

\subparagraph{First statement of the lemma:}
For $n\geq 1$, write $x_n=\lambda_n u_n+v_n$, with $\lambda_n\in\R$ and $v_n\in u_n^\perp$ and denote by $K_n=K+\lambda_n u_n$. Then, by the triangle inequality,
\begin{align} 
	& \textsf{Vol}_{k-1}\left(\left((K+x_n)\cap u_n^\perp\right)\triangle (K\cap u_n^\perp)\right) \leq  \nonumber \\ 
	& \hspace{20mm} \textsf{Vol}_{k-1}\left(\left((K_n+v_n)\cap u_n^\perp\right)\triangle (K_n\cap u_n^\perp)\right) \nonumber \\ 
	\label{TriangleIneq} & \hspace{25mm} +\textsf{Vol}_{k-1}\left(\left((K+\lambda_n u_n)\cap u_n^\perp\right)\triangle (K\cap u_n^\perp)\right).
\end{align}
Since $v_n\in u_n^\perp$, the first term on the right hand side of \eqref{TriangleIneq} is equal to 
\begin{equation} \label{EqnVolumes}
	\textsf{Vol}_{k-1}\left(\left((K_n\cap u_n^\perp)+v_n\right)\triangle (K_n\cap u_n^\perp)\right).
\end{equation}
It is easy to see that $\DS d_{\textsf{H}}\left((K_n\cap u_n^\perp)+v_n, K_n\cap u_n^\perp\right)\leq |v_n|$, which is less than one if $n$ is large enough. Hence, using the same argument as in the proof of Lemma 1 in \citep{Brunel2013}, there is a positive constant $C$ that does not depend on $n$ such that \eqref{EqnVolumes} is bounded from above by $\DS Cd_{\textsf{H}}\left((K_n\cap u_n^\perp)+v_n, K_n\cap u_n^\perp\right)$. Therefore, the first term of the right hand side of \eqref{TriangleIneq} goes to zero as $n$ goes to infinity. 
Let $n\geq 1$ be large enough so $\lambda_n< m$. Set $\alpha_n=p_K(u_n)$ and $\beta_n=p_K(-u_n)$. Suppose that $\lambda_n\geq 0$ (the case $\lambda_n< 0$ would be handled similarly). Then, by convexity of $K$, 
\begin{equation} \label{InclusionTruc0000} 
	\frac{\alpha_n}{\alpha_n+\lambda_n}\left(K+\lambda_n u_n\right)\subseteq K \subseteq \frac{\beta_n}{\beta_n-\lambda_n}\left(K+\lambda_n u_n\right).
\end{equation}
Since $0\in K$, it is true that for all $\lambda\in \R$ and $v\in\mathcal S^{k-1}$, $(\lambda K)\cap v^{\perp}=\lambda(K\cap v^\perp)$. Using this fact together with \eqref{InclusionTruc0000} yields
\begin{align}
	& \left((K+\lambda_n u_n)\cap u_n^\perp\right)\triangle (K\cap u_n^\perp) \subseteq \nonumber \\
	& \hspace{20mm} \left(\left(\frac{\alpha_n+\lambda_n}{\alpha_n}(K\cap u_n^\perp)\right)\setminus (K\cap u_n^\perp)\right) \nonumber \\
	\label{InclusionTruc1} & \hspace{25mm} \cup \left((K\cap u_n^\perp)\setminus \left(\frac{\beta_n-\lambda_n}{\beta_n}(K\cap u_n^\perp)\right)\right).
\end{align}
Since $0\in K$, the volume of the set in the right hand side of \eqref{InclusionTruc1} is bounded from above by 
$$\left(\left(\frac{\alpha_n+\lambda_n}{\alpha_n}\right)^{k-1}-1+\left(\frac{\beta_n-\lambda_n}{\beta_n}\right)^{k-1}-1\right)\textsf{Vol}_{k-1}(K\cap u_n^\perp),$$
which goes to zero as $n$ goes to infinity, since $K$ is bounded and $\alpha_n$ and $\beta_n$ are bounded away from zero (they are at not smaller than $m$). This ends the proof of the first statement of the lemma.

\subparagraph{Second statement of the lemma:}
Since the convex set $K$ is bounded, its gauge function is Lipschitz on $\mathcal S^{k-1}$ and it is bounded away from zero on $\mathcal S^{k-1}$. Hence, its (multiplicative) inverse $p_K$ is also Lipschitz on $\mathcal S^{k-1}$. Let $L$ be the corresponding Lipschitz constant. Let $t_n=\|U_n-I_k\|$, where $I_k$ is the identity map in $\R^k$ and we define the norm of any linear map $A:\R^k \to \R^k$ by $\DS \|A\|=\max_{v\in\mathcal S^{k-1}}|A(v)|$. Then, since $U_n$ converges to the identity, $t_n$ goes to zero as $n$ goes to infinity. Define $\DS c_n=\frac{m}{m+Lt_n}$. Note that $0\leq c_n\leq 1$. Then, let us show that for all $n\geq 1$,
\begin{equation} \label{InclusionTruc2}
	c_n K \subseteq U_n(K).
\end{equation}
Let $x\in K$ and set $y=U_n^{-1}(c_n x)$. If $x=0$, then $y=0$ yielding $y\in K$ by assumption, which proves \eqref{InclusionTruc2}. If $x\neq 0$, then $y\neq 0$ and let $v=y/|y|$. In order to prove that $y\in K$, it is enough to show that 
\begin{equation} \label{InclusionTruc3}
	|y|\leq p_K(v).
\end{equation}
Since $U_n$ is an isometry, $|y|=c_n|x|$ and since $x\in K$, $|x|\leq p_K(x/|x|)$. Therefore,
\begin{align*}
	|y| & = c_n|x| \leq c_n p_K(x/|x|) = c_n p_K(U_n(v)) \leq c_n p_K(v) + c_n L|U_n(v)-v|  \\
	& \leq c_n p_K(v) + c_n L t_n = p_K(v) + c_n L t_n - (1-c_n)p_K(v)  \leq p_K(v) + c_n t_n L - (1-c_n)m \\
	& = p_K(v),
\end{align*}
by definition of $c_n$. This proves \eqref{InclusionTruc3} and hence, \eqref{InclusionTruc2}. As a consequence, since $0\in K$, $\DS c_n (K\cap u_n^\perp) = (c_n K)\cap u_n^\perp \subseteq U_n(K)\cap u_n^\perp$, yielding 
\begin{align*}
	\textsf{Vol}_{k-1}\left((K\cap u_n^\perp)\setminus\left(U_n(K)\cap u_n^\perp\right) \right) & \leq \textsf{Vol}_{k-1}\left((K\cap u_n^\perp)\setminus \left(c_n^{-1}(K\cap u_n^\perp)\right)\right) \\
	& \leq (1-c_n^{-(k-1)})\textsf{Vol}_{k-1}(K\cap u_n^\perp),
\end{align*}
which goes to zero as $n\to\infty$, since $K$ is bounded and $c_n\to 1$. In a similar fashion, we prove that $\displaystyle{\textsf{Vol}_{k-1}\left(\left(U_n(K)\cap u_n^\perp\right)\setminus(K\cap u_n^\perp) \right)}$ also goes to zero as $n\to\infty$, which ends the proof of the second statement of the lemma. \hfill \textsquare

%
%
%
%
%
%
%
%

\begin{lemma} \label{LemmaBarB}
Let $M$ be a positive integer and let $U_1,\ldots,U_M$ be i.i.d. uniform random variables on $\mathcal S^{d-1}$. Let $\delta\in (0,1]$ and let $\mathcal C$ be the event satisfied when the collection $\{U_1,\ldots,U_M\}$ is a $\delta$-net of the sphere (see Definition \ref{DefNets}). Then, the complement $\mathcal C^{\complement}$ of $\mathcal C$ satisfies $\DS \PP[\mathcal C^{\complement}] \leq \#\mathcal N \left(1-\left(\frac{\delta}{4}\right)^{(d-1)/2}\right)^M  \leq 6^d\exp\left(-\frac{M\delta^{d-1}}{2d8^{(d-1)/2}}+d\log\left(\frac{1}{\delta}\right)\right)$.
\end{lemma}

\paragraph{Proof:} Let $\mathcal N$ be a $(\delta/2)$-net of $\mathcal S^{d-1}$. By a simple volume argument, it is possible to choose $\mathcal N$ satisfying $\#\mathcal N\leq (6/\delta)^d$, which we assume in the sequel. If $\mathcal C$ is not satisfied, there exists $u\in\mathcal S^{d-1}$ for which $|u-U_j|>\delta$, for all $j=1,\ldots,M$. Hence, if $v\in\mathcal N$ is such that $|u-v|\leq \delta/2$, one has, for all $j=1,\ldots,M$, by the triangle inequality, $\DS |v-U_j| \geq |u-U_j|-|u-v| \geq \delta-\delta/2 \geq \delta/2$. Therefore, using the union bound and mutual independence of the $U_j$'s,
\begin{equation}\label{label1}
	\PP[\mathcal C^{\complement}] \leq \PP\left[\exists v\in\mathcal N, |v-U_j|>\frac{\delta}{2}, \forall j=1,\ldots,M\right] \leq \sum_{v\in\mathcal N} \PP\left[|v-U_1|>\frac{\delta}{2}\right]^M.
\end{equation}
For any $v\in\mathcal S^{d-1}$, $\displaystyle{\PP\left[|v-U_1|\leq\frac{\delta}{2}\right]}$ is the ratio of the surface area of a spherical cap of the unit sphere and the total surface area of the unit sphere. The height of this cap is $h=\delta^2/8<1$. Then,
\begin{equation} \label{label3}
	\PP\left[|v-U_1|\leq\frac{\delta}{2}\right] = \frac{1}{2}I_{2h-h^2}\left(\frac{d-1}{2},\frac{1}{2}\right),
\end{equation}
where $\DS I_x(a,b)=\frac{\int_0^x t^{a-1}(1-t)^{b-1}\diff t}{\int_0^1 t^{a-1}(1-t)^{b-1}\diff t}$, for $x\in [0,1]$ and $a,b>0$. If $b\leq 1$, one has $\DS \int_0^x t^{a-1}(1-t)^{b-1}\diff t  \geq \int_0^x t^{a-1} \diff t = \frac{x^a}{a}$ and $\DS a\int_0^1 t^{a-1}(1-t)^{b-1}\diff t = (a+b)\int_0^1 t^{a}(1-t)^{b-1} \diff t \leq (a+b)\int_0^1 (1-t)^{b-1} \diff t = \frac{a+b}{b}$. Hence, $\displaystyle{I_x(a,b)\geq \frac{b}{a+b} x^a}$ and \eqref{label3} yields, with $x=2h-h^2$, $a=\frac{d-1}{2}$ and $b=1/2$, that $\DS \PP\left[|v-U_1|\leq\frac{\delta}{2}\right] \geq \frac{1}{2d}(2h-h^2)^{(d-1)/2}$. Since $h<1$, $2h-h^2\geq h=\delta^2/8$, hence,
\begin{equation} \label{label5}
	\PP\left[|v-U_1|\leq\frac{\delta}{2}\right] \geq \frac{\delta^{d-1}}{2d8^{(d-1)/2}}.
\end{equation}
Together with \eqref{label5}, \eqref{label1} implies
\begin{equation*}
	\PP[\mathcal C^{\complement}] \leq \#\mathcal N \left(1-\left(\frac{\delta}{4}\right)^{(d-1)/2}\right)^M  \leq 6^d\exp\left(-\frac{M\delta^{d-1}}{2d8^{(d-1)/2}}+d\log\left(\frac{1}{\delta}\right)\right),
\end{equation*}
which ends the proof of Lemma \ref{LemmaBarB}. \hfill \textsquare

\subsection{Preliminary lemmas for empirical and population quantiles}

\begin{lemma} \label{ContTrueQuant}

Let $\mu$ satisfy Assumption \ref{AssA}. Then, the map $u\in\mathcal S^{d-1}\mapsto q_u^\sharp$ is continuous.

\end{lemma}

\paragraph{Proof: }

For notation's sake, we write $q_u$ instead of $q_u^\sharp$ in the sequel of the proof. 
\subparagraph{Step 1:} Denote by $\Phi(u,t)=\PP[\langle u,X\rangle\leq t], u\in\mathcal S^{d-1}, t\in\R$. We first show that $\Phi$ is continuous $\displaystyle{A=\left\{(u,t)\in\mathcal S^{d-1}\times\R:q_u-\varepsilon< t< q_u+\varepsilon\right\}}$.

Let $(u,t)\in A$ and $(u_p,t_p)_{p\geq 1}$ be a sequence in $A$ that converges to $(u,t)$ as $p$ goes to infinity. Let $\eta$ be an arbitrary positive number. We show that if $p$ is large enough, then $\DS \left|\Phi(u_p,t_p)-\Phi(u,t)\right|\leq 2\eta$,
which will prove our statement. First, note that 
$\DS \left|\Phi(u_p,t_p)-\Phi(u,t)\right|\leq \mu\left(H_{u,t}\triangle H_{u_p,t_p}\right)$.
Let $R>0$ satisfy $\DS\PP[|X|>R]\leq \eta$. Then, 
\begin{align*}
\mu\left(H_{u,t}\triangle H_{u_p,t_p}\right) & \leq \mu\left(B(0,R)\cap(H_{u,t}\triangle H_{u_p,t_p})\right)+\mu(R^d\setminus B(0,R)) \\
& \leq \mu\left(B(0,R)\cap(H_{u,t}\triangle H_{u_p,t_p})\right)+\eta.
\end{align*}
It is easy to check that
$$B(0,R)\cap(H_{u,t}\triangle H_{u_p,t_p})\subseteq \left(H_{u,t_p+R|u_p-u|}\setminus H_{u,t}\right)\cup \left(H_{u,t}\setminus H_{u,t_p-R|u_p-u|}\right),$$
which entails
\begin{align}
& \mu\left(B(0,R)\cap(H_{u,t} \triangle H_{u_p,t_p})\right) \nonumber \\
\label{eqn12345} & \hspace{15mm}\leq \left|F_u\left(t_p+R|u_p-u|\right)-F_u(t)\right|+\left|F_u(t)-F_u\left(t_p-R|u_p-u|\right)\right|.
\end{align}
Since $\displaystyle{(u_p,t_p)\xrightarrow[p\to\infty]{} (u,t)}$ and $q_u-\varepsilon< t< q_u+\varepsilon$, one has $\DS q_u-\varepsilon\leq t_p-R|u_p-u|\leq t_p+R|u_p-u|\leq q_u+\varepsilon$ for all large enough $p$. Hence, since $F_u$ is continuous on $[q_u-\varepsilon,q_u+\varepsilon]$, \eqref{eqn12345} implies that $\DS \mu\left(B(0,R)\cap(H_{u,t} \triangle H_{u_p,t_p})\right)\leq\eta$ if $p$ is large enough, which ends the the proof of the continuity of $\Phi$ on $A$. 

\subparagraph{Step 2:} Let $u\in\mathcal S^{d-1}$ and $(u_p)_{p\geq 1}$ be a sequence of unit vectors converging to $u$ as $p$ goes to infinity. Let us show that $q_{u_p}$ converges to $q_u$. If this was not the case, there would be a positive number $\eta$ and an increasing sequence of positive integers $(p_k)_{k\geq 1}$ satisfying $\DS |q_{u_{p_k}}-q_u|\geq\eta, \forall k\geq 1$. Let us assume that $q_{u_{p_k}}\geq q_u+\eta$ for an infinite number of indices $k\geq 1$. The case when $q_{u_{p_k}}\leq q_u-\eta$ for an infinite number of indices $k\geq 1$ would be handled similarly. For the sake of notation, we renumber the sequence and assume that for $k\geq 1$, $q_{u_k}\geq q_u+\eta$. Without loss of generality, assume that $\eta<\varepsilon$. Hence, for all $k\geq 1$,
\begin{align}
1-\alpha=F_{u_k}(q_{u_k}) & \geq F_{u_k}(q_u+\eta) \nonumber \\
& = F_u(q_u+\eta)+ \Phi(u_k,q_u+\eta)-\Phi(u,q_u+\eta) \nonumber \\
& \geq F_u(q_u) + L\eta+\Phi(u_k,q_u+\eta)-\Phi(u,q_u+\eta) \nonumber \\
\label{eqn13472}& = 1-\alpha + L\eta+\Phi(u_k,q_u+\eta)-\Phi(u,q_u+\eta).
\end{align}
The fact that $F_v(q_v)=1-\alpha, \forall v\in\mathcal S^{d-1}$, is a consequence of the continuity and strict monotony of $F_v$ in a neighborhood of $q_v$, for all $v\in\mathcal S^{d-1}$. Since $\eta<\varepsilon$, $(u,q_u+\eta)\in A$, so by the first part of the proof, $\DS \Phi(u_k,q_u+\eta)-\Phi(u,q_u+\eta) \xrightarrow[k\to\infty]{}0$. Thus, by letting $k$ grow to infinity in \eqref{eqn13472}, we get that $L\eta\leq 0$, which is a contradiction. Hence, we have proved that $\displaystyle{q_{u_p}\xrightarrow[p\to\infty] {}q_u}$, which ends the proof. \hfill \textsquare

\begin{lemma} \label{Contphi}

	Let $\mu$ be a probability measure on $\R^d$ that satisfies either Assumption \ref{Ass11} or \ref{Ass12} and let $K$ be its support. For $u\in\mathcal S^{d-1}$, let $f_u$ and $F_u$ be, respectively, the density and the cumulative distribution function of $\langle u,X\rangle$, where $X$ is a random variable with distribution $\mu$. Let $A=\{(u,t)\in \mathcal S^{d-1}\times \R: (tu+u^\perp)\cap K \neq\emptyset\}$. Define 
$\DS \phi(u,t)=f_u(t)$ and $\DS \Phi(u,t)=F_u(t)$, for all $(u,t)\in A$. Then,
\begin{itemize}
	\item $\phi$ and $\Phi$ are continuous on $\overset{\circ}{A}$;
	\item $\forall (u,t)\in\overset{\circ}{A}$, $\phi(u,t)>0$ and $0<\Phi(u,t)<1$.
\end{itemize}
\end{lemma}

\paragraph{Proof:}

Note $\DS \Phi(u,t)=\int_{-\infty}^t\phi(u,s)\diff s$, for all $(u,t)\in A$, where we set $\phi(u,s)$ to zero if $(u,s)\notin A$. Hence, by dominated convergence, continuity of $\phi$ will automatically yield that of $\Phi$. Let $(u,t)\in A$ and consider an arbitrary sequence $(u_n,t_n)_{n\geq 1}$ of elements of $A$ that converges to $(u,t)$.

\subparagraph{Let $\mu$ satisfy Assumption \ref{Ass11}.}

In this case, the second statement is trivial since $f$ is continuous and positive everywhere. Hence, we only prove the first statement. Let $\epsilon>0$ and $R>0$. For all $n\geq 1$,
\begin{align}
	\left|\phi(u_n,t_n)-\phi(u,t)\right| & = \left|\int_{u_n^\perp}f(t_nu_n+v)\diff v - \int_{u^\perp}f(tu+v)\diff v\right| \nonumber \\
	& \leq \left|\int_{\substack{v\in u_n^\perp: \\ |v|\leq R}}f(t_nu_n+v)\diff v - \int_{\substack{v\in u^\perp: \\ |v|\leq R}}f(tu+v)\diff v\right| \nonumber \\
	& \hspace{10mm} +\int_{\substack{v\in u_n^\perp: \\ |v|> R}}f(t_nu_n+v)\diff v + \int_{\substack{v\in u^\perp: \\ |v|> R}}f(tu+v)\diff v \nonumber \\
	& \leq \left|\int_{\substack{v\in u_n^\perp: \\ |v|\leq R}}f(t_nu_n+v)\diff v - \int_{\substack{v\in u^\perp: \\ |v|\leq R}}f(tu+v)\diff v\right| \nonumber \\ 
	\label{Cont1} & \hspace{10mm} + C\int_{\substack{v\in u_n^\perp: \\ |v|> R}}(1+|v|)^{-\nu}\diff v + \int_{\substack{v\in u^\perp: \\ |v|> R}}(1+|v|)^{-\nu}\diff v.
\end{align}

For $n\geq 1$, let $U_n$ be an isometry in $\R^d$ such that $U_n(u_n)=u$. Then, the first term in \eqref{Cont1} can also be written as $\DS \left|\int_{\substack{v\in u^\perp: \\ |v|\leq R}}\left(f(t_nu_n+U_n^{-1}(v))-f(tu+v)\right)\diff v\right|$, which converges to zero by dominated convergence. Hence, for large $n$, the first term in \eqref{Cont1} is smaller than $\epsilon$.

Using polar coordinates, both the second and third terms in \eqref{Cont1} can be rewritten as $\DS C'\int_R^\infty x^{d-2}(1+x)^{-\nu}\diff x,$
for some positive constant $C'$ that does not depend on $n$ or $R$. Hence, both the second and third terms in \eqref{Cont1} are bounded from above by $C''R^{-(\nu-d+1)}$, for some positive constant $C''$ that does not depend on $R$ or $n$. Hence, if $R$ was chosen large enough, both these terms are smaller than $\epsilon$. Finally, we have proved that $\phi(u_n,t_n)\to\phi(u,t)$, as $n\to\infty$.

\subparagraph{Let $\mu$ satisfy Assumption \ref{Ass12}.}
Let $\epsilon>0$. For $n\geq 1$, write
\begin{align} 
	|\phi(u_n,t_n)-\phi(u,t)| & \leq \int_{u_n^\perp}\left|f(t_nu_n+v)-f(tu+v)\right|\diff v \nonumber \\
	\label{Cont2} & \hspace{15mm} + \left|\int_{u_n^\perp}f(tu+v)\diff v - \int_{u^\perp}f(tu+v)\diff v \right|.
\end{align}
Let $\displaystyle{B_n=\{v\in u_n^\perp:t_nu_n+v\in K\}}$ and $\displaystyle{D_n=\{v\in u_n^\perp:tu+v\in K\}}$. The first integral in \eqref{Cont2} can be decomposed as
\begin{align}
	& \int_{B_n\cap D_n}\left|f(t_nu_n+v)-f(tu+v)\right|\diff v \nonumber \\
	\label{Cont3} & \hspace{10mm} + \int_{B_n\triangle D_n}\left|f(t_nu_n+v)-f(tu+v)\right|\diff v.
\end{align}
Recall that $f$ is uniformly continuous on $K$ and $\textsf{Vol}_{d-1}(B_n\cap D_n)$ is bounded uniformly in $n$, by boundedness of $K$. Hence, if $n$ is large enough, the first integral in \eqref{Cont3} is smaller than $\epsilon$. For the second integral, since $f$ is uniformly continuous on the bounded set $K$ and vanishes everywhere else, it is bounded and the integral is bounded from above by $(\sup_K f)\textsf{Vol}_{d-1}(B_n\triangle D_n)$. The latter converges to zero as $n$ goes to infinity, thanks to Lemma \ref{Covariogram}. Hence, it becomes smaller than $\epsilon$ if $n$ is large enough, so the first term in \eqref{Cont2} is at most $2\varepsilon$ for large values of $n$. 
For $n\geq 1$, let $U_n$ be an isometry in $\R^d$ such that $U_n(u_n)=u$ and such that $U_n$ converges to the identity, as $n$ goes to infinity. Then, the second term in the right hand side of \eqref{Cont2} can be written as 
\begin{equation} \label{Cont2SecondTerm}
	\left|\int_{u^\perp} \left(f(tu+U_n^{-1}(v))-f(tu+v)\right)\diff v\right|.
\end{equation}
Let $K_u=(K-tu)\cap u^{\perp}$ and $K_u^{n}=\left(U_n^{-1}(K-tu)\right)\cap u^{\perp}$. Since the integrand vanishes outside of $K_u\cup K_u^{(n)}$, the integral inside the absolute value in \eqref{Cont2SecondTerm} can be decomposed as the sum of two integrals: One on $K_u\cap K_u^n$ and the other on $K_u\triangle K_u^n$.
Since $U_n$ converges to the identity as $n$ goes to infinity, $U_n^{-1}(v)\to v$ as $n\to\infty$, for all $v\in\R^d$. Since $f$ is uniformly continuous on $K$ and $K$ is bounded, $f$ is bounded. Hence, by dominated convergence, uniform continuity of $f$ on $K$ together with the fact that $\textsf{Vol}_{d-1}(K_u\cap K_u^n)$ is bounded uniformly in $n$ implies that the first term goes to zero as $n\to\infty$. Since $f(x)=0$ for $x\notin K$, $f$ is bounded on $\R^d$. Hence, by Lemma \ref{Covariogram}, the second term goes to zero as $n\to\infty$, since $U_n$ converges to the identity. This ends the proof of the first statement of the lemma. \vspace{3mm}

For the second statement, first note that $K$ needs to have a nonempty interior. Otherwise, since it is convex, it would be included in a hyperplane, i.e., there would exist $u\in\mathcal S^{d-1}$ and $t\in\R$ such that $\langle u,x\rangle=t, \forall x\in K$. Hence, $\langle u,X\rangle=t$ almost surely, which contradicts the fact that $X$ has a density with respect to the Lebesgue measure in $\R^d$. Let $(u,t)\in\overset{\circ}{A}$. By Lemma \ref{LemmaInteriorA}, there exists $\eta>0$ such that both $(t+\eta)u+u^\perp$ and $(t-\eta)u+u^\perp$ intersect $\overset{\circ}{K}$. Hence, by convexity of $\overset{\circ}{K}$, $(su+u^\perp)\cap \overset{\circ}{K}\neq\emptyset$, yielding that the $(d-1)$-dimensional Lebesgue measure of $(su+u^\perp)\cap K$ needs to be positive, for all $s\in [t-\eta/2,t+\eta/2]$. Therefore, $f_u$ is positive on this interval, yielding $\phi(u,t)>0$ and $0<\Phi(u,t)<1$. \hfill \textsquare

\begin{lemma} \label{ContQuant28}

	Let $\mu$ be a probability measure on $\R^d$ that satisfies either Assumption \ref{Ass11} or \ref{Ass12} and let $X$ be a random variable with distribution $\mu$. Denote by $F_u$ the cumulative distribution function of $\langle u,X\rangle$. Let $\beta\in (0,1)$. For $u\in\mathcal S^{d-1}$, let $q_u$ be the $\beta$-quantile of $\langle u,X\rangle$, defined as in Lemma \ref{ContTrueQuant} (with $\beta=1-\alpha$). Then, 
\begin{itemize}
	\item For all $u\in\mathcal S^{d-1}$, $q_u$ is the unique real number $t$ that satisfies $F_u(t)=\beta$;
	\item The map $u\in\mathcal S^{d-1}\mapsto q_u$ is continuous.
\end{itemize}

\end{lemma}

\paragraph{Proof:} Let $K$ be the support of $\mu$ and let $A=\{(u,t)\in \mathcal S^{d-1}\times \R: (tu+u^\perp)\cap K \neq\emptyset\}$.

Let $u\in\mathcal S^{d-1}$. Since $\mu$ is absolutely continuous with respect to the Lebesgue measure, so is the distribution of $\langle u,X\rangle$. Hence, $F_u$ is continuous on $\R$, which yields that $F_u(q_u)=\beta$. In addition, if $f_u$ is the density of $\langle u,X\rangle$, then $f_u$ is positive in a neighborhood of $q_u$. Indeed, since $K$ is convex, the support of $f_u$ is an interval. Since $F_u(q_u)=\beta\in (0,1)$ and $F_u$ is continuous, there is a neighborhood of $q_u$ on which $F_u(t)\in (0,1)$, i.e., there is a neighborhood of $q_u$ that is included in the support of $f_u$. In particular, $F_u$ is strictly increasing on this neighborhood, which shows the uniqueness of $q_u$.

Let $u\in\mathcal S^{d-1}$ and let $(u_n)_{n\geq 1}$ be an arbitrary sequence of unit vectors that converges to $u$. Suppose that $q_{u_n}$ does not converge to $q_u$. Then, there exists $\eta>0$ and a subsequence of $u_n$ (renamed $u_n$ after renumbering) such that $|q_{u_n}-q_u|\geq\eta$, for all $n\geq 1$. Assume that for an infinite number of indices $n$, $q_{u_n}\geq q_u+\eta$. The case when $q_{u_n}\leq q_u-\eta$ for an infinite number of indices $n$ would be handled similarly. Thus, up to renumbering the sequence again, assume that $q_{u_n}\geq q_u+\eta$, for all $n\geq 1$. By a similar argument as in the end of the proof of Lemma \ref{Contphi}, for all $(u,t)\in A$, $(u,t)\in \overset{\circ}{A}$ if and only if $0<F_u(t)<1$. Hence, $(u,q_u)\in \overset{\circ}{A}$. Hence, there exists $\xi>0$ such that $(v,t)\in \overset{\circ}{A}$ for all $v\in\mathcal S^{d-1}$ and $t\in\R$ with $|v-u|\leq \xi$ and $|q_u-t|\leq \xi$. By Lemma \ref{Contphi}, since $\phi$ is continuous and positive on $\overset{\circ}{A}$, there is a positive constant $c$ such that $\phi(v,t)\geq c>0$ for all $(v,t)\in \mathcal S^{d-1}\times \R$ with $|v-u|\leq \xi$ and $|q_u-t|\leq \xi$. Assume that $\xi\leq \eta$, without loss of generality. Then,
\begin{align*}
	\beta & = F_{u_n}(q_{u_n}) = \Phi(u_n,q_{u_n}) \geq \Phi(u_n,q_u+\eta) \geq \Phi(u_n,q_u+\xi) = \Phi(u_n,q_u)+\int_0^\xi \phi(u_n,t)\diff t \\
	& \geq \Phi(u_n,q_u)+c\xi \to \beta + c\xi,
\end{align*}
as $n$ goes to infinity. This is a contradiction, since $\beta+c\xi>\beta$. Hence, $q_{u_n}$ needs to converge to $q_u$ as $n\to\infty$ and Lemma \ref{ContQuant28} is proven. \hfill \textsquare

\begin{lemma}\label{LemmaUnifDevQuantiles}
	Let $\mu$ satisfy Assumption \ref{AssA}. Then, for all $n\geq 1$ and $z\in\R$ with $\frac{10\sqrt{5(d+1)}}{L\sqrt n}\leq z<\varepsilon$,
$$\PP\left[\sup_{u\in\mathcal S^{d-1}}|\hat q_u-q_u^\sharp|\leq z\right]\geq 1-A\exp\left(-L^2z^2n/2+10\sqrt{5(d+1)}Lz\sqrt n\right),$$
where $A=e^{-250(d+1)}$.
\end{lemma}

\paragraph{Proof: }

Let $\mathcal C_0=\{(u,t)\in\mathcal S^{d-1}\times \R: q_u^\sharp-\varepsilon\leq t\leq q_u^\sharp+\varepsilon\}$ and $\widetilde{\mathcal C}_0= \{(u,t)\in\mathcal C_0:u\in\Q^{d-1}\times \Q, t\in\Q\}$. Denote by $\mathcal H_0=\{H_{u,t}:(u,t)\in\mathcal C_0\}$ and $\widetilde{\mathcal H}_0=\{H_{u,t}:(u,t)\in\widetilde{\mathcal C}_0\}$.

\subparagraph{Step 1:}  We first show that 
\begin{equation} \label{Step1UnifDev}
	\sup_{H\in\mathcal H_0}|\mu_n(H)-\mu(H)|=\sup_{H\in\widetilde{\mathcal H}_0}|\mu_n(H)-\mu(H)| \quad \mbox{almost surely.}
\end{equation}  
If $(u,t)\in\mathcal S^{d-1}$, denote by $\hat F_u(t)=\mu_n(H_{u,t})$, i.e., the empirical cumulative distribution function of $\langle u,X\rangle$. Then, 
\begin{align*}
	\sup_{H\in\mathcal H_0} |\mu_n(H)-\mu(H)| & = \sup_{(u,t)\in\mathcal C_0}|\hat F_u(t)-F_u(t)| \\
	 & = \max\left(\sup_{(u,t)\in\mathcal C_0}(\hat F_u(t)-F_u(t)),\sup_{(u,t)\in\mathcal C_0} (F_u(t)-\hat F_u(t)) \right).
\end{align*}
Hence, it suffices to prove that $\DS \sup_{(u,t)\in\mathcal C_0}(\hat F_u(t)-F_u(t))=\sup_{(u,t)\in\widetilde{\mathcal C}_0}(\hat F_u(t)-F_u(t))$ and that $\DS \sup_{(u,t)\in\mathcal C_0} (F_u(t)-\hat F_u(t))=\sup_{(u,t)\in\widetilde{\mathcal C}_0} (F_u(t)-\hat F_u(t))$. The first statement follows from two facts. First, $\widetilde{\mathcal C}_0$ is dense in $\mathcal C_0$. Second, $(u,t)\mapsto \hat F_u(t)$ is lower semicontinuous and $(u,t)\in\mathcal C_0\mapsto \hat F_u(t)$ is continuous on $\mathcal C_0$, as proved in Step 1 of the proof of Lemma \ref{ContTrueQuant}, yielding that $(u,t)\in\mathcal C_0\mapsto \hat F_u(t)-F_u(t)$ is lower semicontinuous on $\mathcal C_0$. For the second statement, note that for all $u\in\mathcal S^{d-1}$, continuity of $F_u$ on the segment $[q_u^\sharp-\varepsilon,q_u^\sharp+\varepsilon]$ implies that $\DS \sup_{q_u^\sharp-\varepsilon<t<q_u^\sharp+\varepsilon}\hat F_u(t)-F_u(t)=\sup_{q_u^\sharp-\varepsilon<t<q_u^\sharp+\varepsilon}\hat G_u(t)-G_u(t)$, where $G_u(t)=\PP[\langle u,X\rangle\geq t$ and $\DS \hat G_u(t)=\frac{1}{n}\sum_{i=1}^n \mathds 1_{\langle u,X_i\rangle\geq t}$. Then, the same argument as above yields the second statement, and proves \eqref{Step1UnifDev}. In particular, the random variable $\DS \sup_{(u,t)\in\mathcal C_0}|\hat F_u(t)-F_u(t)|$ is measurable and the probability term in the statement of the lemma is well defined.

\subparagraph{Step 2: } Let $u\in\mathcal S^{d-1}$. By definition of $\hat q_u$, the following holds for all $t\in\R$, where, as we recall, $H_{-u,-t}$ is the halfspace $H=\{x\in\R^d:\langle u,x\rangle\geq t\}$:
\begin{itemize}
	\item If $t< \hat q_u$, then $\mu_n(H_{-u,-t})\geq \alpha$,
	\item If $t>\hat q_u$, then $\mu_n(H_{-u,-t})<\alpha$.
\end{itemize}

Assume that for some $u\in\mathcal S^{d-1}$, $|\hat q_u-q_u^\sharp|> z$. Then, either $\hat q_u>q_u^\sharp+z$ or $\hat q_u<q_u^\sharp-z$. If $\hat q_u>q_u^\sharp+z$, let $H=H_{-u,-(q_u^\sharp+z)}\in\mathcal H_0$. Then, $\DS \mu_n(H)\geq\alpha$. Hence, by Assumption \ref{AssA}, $\DS \mu(H) = \PP[\langle u,X\rangle\geq q_u^\sharp+z] = 1-F_u(q_u^\sharp+z) \leq 1-F_u(q_u^\sharp)-Lz = \alpha-Lz$, yielding that $\DS \mu_n(H)-\mu(H)\geq Lz$. If $\hat q_u<q_u^\sharp-z$ a similar reasoning yields $\DS |\mu(H)-\mu_n(H)\geq Lz$ for $H=H_{-u,-(q_u^\sharp-z)}\in\mathcal H_0$. Hence, using \eqref{Step1UnifDev}, it follows that 
\begin{equation} \label{EmpProc1}
	\PP\left[\sup_{u\in\mathcal S^{d-1}}|\hat q_u-q_u^\sharp|>z\right]\leq \PP\left[\sup_{H\in\widetilde{\mathcal H}_0}|\mu_n(H)-\mu(H)|\geq Lz\right].
\end{equation}
Now, denote by $\DS S=\sup_{H\in\widetilde{\mathcal H}_0}|\mu_n(H)-\mu(H)|$. Since $\widetilde{\mathcal H}_0\subseteq \mathcal H$, it has Vapnik-Chervonenkis dimension at most $d+1$. Moreover, it is a countable class of sets, so Proposition 3.1 in \cite{Baraud2016} yields $\E[S]\leq \frac{10\sqrt{5(d+1)}}{\sqrt n}$. Therefore, by Theorem 2.5 in \cite{Koltchinskii2008}, if $Lz\geq \frac{10\sqrt{5(d+1)}}{\sqrt n}$,
\begin{align}
	\PP[S\geq Lz] & \leq \PP\left[S-\E[S]\geq Lz-\frac{10\sqrt{5(d+1)}}{\sqrt n}\right] \nonumber \\
	\label{EmpProcTh} & \leq A\exp\left(-L^2z^2n/2+10\sqrt{5(d+1)}Lz\sqrt n\right),
\end{align}
where $A=e^{-250(d+1)}$. Lemma \ref{LemmaUnifDevQuantiles} follows from \eqref{EmpProc1} and \eqref{EmpProcTh}. \hfill \textsquare

\begin{lemma}\label{LemmaOrderFunctions}

Let $f_1,\ldots,f_n$ be $n$ real valued continuous functions defined on a topological space $E$ and $k\in\{1,\ldots,n\}$. For $x\in E$, denote by $f_{(k)}(x)$ the $k$-th smaller number in the list $f_1(x),\ldots,f_n(x)$. Then, $f_{(k)}$ is continuous.

\end{lemma}

\paragraph{Proof:} Write $\DS f_{(k)}(x)=\min_{J\in\mathcal P_k}\max_{j\in J} f_j(x)$, where $\mathcal P_k$ is the collection of all subsets of $\{1,\ldots,n\}$ of size $k$. Continuity of $f_{(k)}$ follows from continuity of the maximum and minimum of finitely many continuous functions. \hfill \textsquare

\subsection{Proofs of the main theorems}

\paragraph{Proof of Lemma \ref{ThmEqualitySets}: } 

Let us first show that $G^{\flat}_{\textsf{MQ}}=G_{\textsf{FB}}$. Let $x\in G^{\flat}_{\textsf{MQ}}$ and $H\in\mathcal H$ satisfying $\mu(H)\geq 1-\alpha$. Write $H=H_{u,t}$, for some $u\in\mathcal S^{d-1}$ and $t\in\R$. Then, $\DS \mu(H)=\PP[\langle u,X\rangle \leq t]\geq 1-\alpha$, which yields $t\geq q_u^{\flat}$. Since $x\in G^{\flat}_{\textsf{MQ}}$, $\langle u,x\rangle\leq q_u^\flat$ and, hence, $x\in\ H$. Therefore, $G^{\flat}_{\textsf{MQ}}\subseteq G_{\textsf{FB}}$. Now, let $x\in G_{\textsf{FB}}$ and $u\in\mathcal S^{d-1}$. Let $H=H_{u,q^{\flat}_u}$. By definition of $q^{\flat}_u$ and since $F_u$ is right continuous, $\mu(H)=F_u(q^{\flat}_u)\geq 1-\alpha$, so $x\in H$. Hence, $x\in G_{\textsf{MQ}}^\flat$ and thus, $G_{\textsf{FB}}\subseteq G^{\flat}_{\textsf{MQ}}$. This ends the proof of the equality $G^{\flat}_{\textsf{MQ}}=G_{\textsf{FB}}$. \newline
Inclusion $G^{\flat}_{\textsf{MQ}}\subseteq G^{\sharp}_{\textsf{MQ}}$ follows from the inequalities $q^{\flat}_u\leq q^{\sharp}_u$, for all $u\in\mathcal S^{d-1}$. \newline
Now, let us prove that $G^{\sharp}_{\textsf{MQ}}=G_\mu$. For $x\in G^{\sharp}_{\textsf{MQ}}$, we show that $D_{\mu}(x)\geq\alpha$, i.e., that any closed halfspace $H$ containing $x$ needs to satisfy $\mu(H)\geq\alpha$. Let $H$ be such a halfspace and write $H=H_{u,t}$ for some $u\in\mathcal S^{d-1}$ and $t\in\R$. Then, $\langle u,x\rangle \leq t$,so $\langle -u,x\rangle \geq -t$. Since $x\in G^{\sharp}_{\textsf{MQ}}$, $\langle -u,x\rangle\leq q^{\sharp}_{-u}$, hence, $-t\leq q^{\sharp}_{-u}$. Therefore,
\begin{align*}
\mu(H) & = \PP[\langle u,X\rangle\leq t] = 1-\PP[\langle u,X\rangle > t] = 1-\PP[\langle -u,X\rangle < -t] \\ & \geq  1-\PP[\langle -u,X\rangle < q^{\sharp}_{-u}] \geq 1-(1-\alpha) = \alpha.
\end{align*}
Thus, $x\in G_\mu$, and hence, $G^{\sharp}_{\textsf{MQ}}\subseteq G_\mu$. Now, let $x\in G_\mu$ and $u\in\mathcal S^{d-1}$. Since $x\in H_{-u,\langle -u,x\rangle}$ and $D_{\mu}(x)\geq\alpha$, $\mu(H_{-u,\langle -u,x\rangle})\geq\alpha$, i.e., $\PP[\langle -u,X\rangle\leq \langle -u,x\rangle] \geq\alpha$. Hence, 
$\DS \PP[\langle u,X\rangle < \langle u,x\rangle] \leq 1-\alpha$,which, by definition of $q^{\sharp}_u$, implies that $\langle u,x\rangle\leq q^{\sharp}_u$. So, $x\in G^{\sharp}_{\textsf{MQ}}$. Therefore, $G^{\sharp}_{\textsf{MQ}}=G_\mu$. \hfill \textsquare

\paragraph{Proof of Proposition \ref{LemmaPolyhedralRep}}

\begin{itemize}
	\item \textit{(i)} $\Rightarrow$ \textit{(ii)}: Assume that all the constraints are active and let $u\in\mathcal S^{d-1}$. First, by definition of the support function, $h_G(u)\leq t_u$. Second, since the constraint corresponding to $u$ is active, there exists $x^*\in G$ such that $\langle u,x^*\rangle=t_u$, yielding $t_u\leq h_G(u)$, hence, $t_u=h_G(u)$.
	\item \textit{(ii)} $\Rightarrow$ \textit{(i)}: Let $u\in\mathcal S^{d-1}$. By Lemma \ref{Lemmacb}, $G$ is compact, yielding the existence of $x^*\in G$ satisfying $h_G(u)=\langle u,x^*\rangle$. Hence, the constraint corresponding to $u$ is active.
	\item \textit{(ii)} $\Rightarrow$ \textit{(iii)} is a direct consequence of the sublinearity of support functions.
	\item \textit{(iii)} $\Rightarrow$ \textit{(ii)}: Assume that the family $(t_u)_{u\in\R^d}$ is sublinear and let $u_0\in\mathcal S^{d-1}$. Since $u\in\R^d \mapsto t_u$ is sublinear and positively homogeneous, it is convex. Hence, it is continuous on the interior of its domain, here, $\R^d$. Since $\overset{\circ}{G}\neq\emptyset$, Lemma \ref{LemmaActiveCone} yields the existence of $x^*\in G$ satisfying $h_G(u_0)=\langle u_0,x^*\rangle$ and of $u_1,\ldots,u_d\in\mathcal S^{d-1}$, $\lambda_1,\ldots,\lambda_d\geq 0$ satisfying $u_0=\sum_{i=1}^d \lambda_i u_i$ and, for $i=1,\ldots,d$, $\langle u_i,x^*\rangle = t_{u_i}$. Hence, $\DS
	h_G(u_0) = \langle u_0,x^*\rangle = \sum_{i=1}^d \lambda_i \langle u_i,x^*\rangle = \sum_{i=1}^d \lambda_i t_{u_i} \geq t_{u_0}$, by positive homogeneity and sublinearity of $v\mapsto t_v$. Since, in addition, $h_G(u_0)\leq t_{u_0}$ by definition of the support function, $h_G(u_0)\leq t_{u_0}$.
\end{itemize}
\hfill \textsquare

\paragraph{Proof of Theorem \ref{EmptyLevelSets}}

Let $\alpha$ be greater than $1/2$ and assume that $G_\mu$ is nonempty. Let $x\in G_\mu$: We prove that $\mu(\{x\})>0$. 

Let $E$ be an affine hyperplane passing through $x$. Let $H_1$ and $H_2$ be the two distinct halfspaces whose common boundary is $E$. Since $x\in G_\mu$, $D_\mu(x)\geq \alpha$. In particular, since both $H_1$ and $H_2$ contain $x$, $\mu(H_j)\geq \alpha, j=1,2$. Hence, $\DS 1\geq \mu(H_1\cup H_2)=\mu(H_1)+\mu(H_2)-\mu(E)\geq 2\alpha-\mu(E)$, which implies that $\mu(E)\geq 2\alpha-1$.

Let $k\in\{1,\ldots,d\}$. Assume it is known that any affine subspace $E$ of dimension $k$, containing $x$, satisfies $\mu(E)\geq 2\alpha-1$.  Let $F$ be an affine subspace of dimension $k-1$, containing $x$. Let $G$ be the linear subspace of vectors that are orthogonal to $F$. Let $p\geq 2$ be an integer and let $u_1,\ldots,u_p$ be unit vectors in $G$, such that no two of them are collinear. For $i=1,\ldots,p$, set $E_i=F+\R u_i=\{f+\lambda u_i:f\in F, \lambda\in\R\}$. Then, for all $I\subseteq \{1,\ldots,p\}$ with $\#I\geq 2$, $\displaystyle{\bigcap_{i\in I}E_i}=F$ and as a consequence of the inclusion-exclusion principle,
\begin{equation*}
	1 \geq \mu\left(\bigcup_{i=1}^p E_i\right) = \sum_{i=1}^p \mu(E_i) - \sum_{j=2}^p (-1)^j {p \choose j}\mu(F) \geq p(2\alpha-1) - (p-1)\mu(F),
\end{equation*}
yielding $\DS \mu(F)\geq \frac{p}{p-1} (2\alpha-1) -\frac{1}{p-1}$. Since $p$ is an arbitrary integer, we can let it go to infinity and we get $\mu(F)\geq 2\alpha-1$.

By induction, this proves that $\mu(\{x\})\geq 2\alpha-1>0$ and this must hold for all $x\in G_\mu$. Since $G_\mu$ is convex, it cannot contain more than one point. Indeed, if $x,y\in G_\mu$, then $[x,y]\subseteq G_\mu$, yielding $\mu(\{z\})\geq 2\alpha-1$, for all $z\in [x,y]$. Hence, if $x\neq y$, then $\mu([x,y])=\infty$, which is impossible. \hfill \textsquare

\paragraph{Proof of Theorem \ref{MainTheorem}}

Without loss of generality, let us assume that $a=0$ in Assumption \ref{AssB}: translating the measure $\mu$ and the sample points does not affect the Haussdorf distance between $G_{\mu}$ and $\hat G$. For the sake of notation, we write $q_u=q_u^{\flat}=q_u^{\sharp}$ for all $u\in\mathcal S^{d-1}$.

Let $z\in [10\sqrt{5(d+1)}/(L\sqrt n),\varepsilon)$ and let the event $\DS \mathcal A=\{|\hat q_u-q_u|\leq z, \forall u\in\mathcal S^{d-1}\}$ hold. Since $B'(0,r)\subseteq G_\mu$, it is true that $q_u\geq r, \forall u\in\mathcal S^{d-1}$. Hence, for all $u\in\mathcal S^{d-1}$, $\hat q_u\geq q_u-z\geq r-\varepsilon>0$, yielding that $B'(0,r-\varepsilon)\subseteq \hat G$, hence, that $\hat G$ has a nonempty interior. So does $G_\mu$, since it contains $B'(0,r)$.


By Lemmas \ref{ContTrueQuant} and \ref{LemmaOrderFunctions}, the maps $u\mapsto q_u$ and $u\mapsto \hat q_u$ are continuous. Indeed, $\hat q_u$ is the $\lceil n(1-\alpha)+1\rceil$-th order function of $\langle u,X_1\rangle, \ldots, \langle u,X_n\rangle$. Note that the map $\DS t\in [0,1)\mapsto \frac{1+t}{1-t}$ is nondecreasing. Thus, by Lemma \ref{Lemma_q_to_h},	$\DS d_{\textsf H}(\hat G,G_\mu)\leq \frac{z R}{r}\frac{1+z/r}{1-z/r} \leq Cz$, where $\DS C=\frac{R}{r}\frac{1+\varepsilon/r}{1-\varepsilon/r}.$ Hence, if $\mathcal A^{{\complement}}$ stands for the complement of the event $\mathcal A$, then
\begin{equation} \label{ProofMainTh1}
	\PP[d_{\textsf H}(\hat G,G_\mu)>Cz]\leq \PP\left[\mathcal A^{{\complement}}\right].
\end{equation}
Write $z=x/\sqrt n$, for some real number $x$ satisfying $\DS \frac{10\sqrt{5(d+1)}}{L}\leq x<\varepsilon\sqrt n$. By Lemma \ref{LemmaUnifDevQuantiles},
\begin{equation} \label{ProofMainTh2}
\PP\left[\mathcal A^{{\complement}}\right]\leq A\exp\left(-L^2x^2/2+10\sqrt{5(d+1)}Lx\right).
\end{equation}
The desired result is a consequence of \eqref{ProofMainTh1} and \eqref{ProofMainTh2}. \hfill \textsquare

\paragraph{Proof of Corollary \ref{CorollFubini}}

The proof is based on a consequence of Funini's theorem which ensures that if $Z$ is a nonnegative random variable, then 
\begin{equation} \label{Fubini_k}
	\E[Z^k]=k\int_0^\infty t^{k-1}\PP\left[Z>t\right]\diff t,
\end{equation}
for all positive number $k$.

First, note that for all $k>0$,
\begin{equation} \label{ProofMainCor0}
	\E\left[d_{\textsf H}(\hat G^*,G_\mu)^k\right] = \E\left[d_{\textsf H}(\hat G^*,G_\mu)\mathds 1_{\hat G\neq\emptyset}\right]+\E\left[d_{\textsf H}(\hat G^*,G_\mu)^k\mathds 1_{\hat G=\emptyset}\right],
\end{equation}
where $\mathds 1$ stands for the indicator function. 

By definition of $\hat G^*$, the second term in the right hand side of \eqref{ProofMainCor0} is equal to $\DS d_{\textsf H}(\{0\},G_\mu)^k\PP\left[\hat G=\emptyset\right]$. First, it is clear that $\DS d_{\textsf H}(\{0\},G_\mu)\leq |a|+R\leq\tau+R$. Second, as we saw in the proof of Theorem \ref{MainTheorem}, $\DS \PP\left[\hat G=\emptyset\right]\leq \PP\left[\mathcal A^{\complement}\right]$ where we set $z=\varepsilon/2$. Hence, by \eqref{ProofMainTh2}, 
\begin{equation} \label{ProofMainCor1}
	\E\left[d_{\textsf H}(\hat G^*,G_\mu)^k\mathds 1_{\hat G=\emptyset}\right] = O\left(n^{-k/2}\right),
\end{equation}
with multiplicative constants that depend on $d,\varepsilon,R,L$ and $\tau$ only.

For the first term of \eqref{ProofMainCor0}, note that if $\hat G\neq\emptyset$, then, since $\hat G^*\subseteq B'(0,\log n)$ and $G_\mu\subseteq B'(a,R)$, $\DS d_{\textsf H}(\hat G^*,G_\mu)\leq |a|+\log n+R\leq \tau+\log n+R$. Denote by $B=\tau+\log n+R$. Then, if we set $Z=d_{\textsf H}(\hat G,G_\mu)$,
\begin{equation} \label{ProofMainCor2}
	\E\left[d_{\textsf H}(\hat G^*,G_\mu)^k\mathds 1_{\hat G\neq \emptyset}\right] \leq \E\left[Z^k\mathds 1_{Z\leq B}\right].
\end{equation}
In the following, we set $k=1$. General values of $k$ would be handled similarly, using \eqref{Fubini_k}.
Using \eqref{ProofMainCor2} and \eqref{Fubini_k} with $k=1$,
\begin{equation*}
	\E[d_{\textsf{H}}(\hat G^*,G_\mu)\mathds 1_{\hat G\neq\emptyset}] \leq \int_0^B \PP\left[d_{\textsf{H}}(\hat G,G_{\mu})>t\right]\diff t.
\end{equation*}
Split the integral in three integrals. First, from $0$ to $\DS \frac{10C\sqrt{5(d+1)}}{L\sqrt n}$, where we bound the integrand by $1$. Second, from $\DS \frac{10C\sqrt{5(d+1)}}{L\sqrt n}$ to $\varepsilon$, where we use the bound provided by Theorem \ref{MainTheorem}. Third, in the remaining interval, where, using monotonicity, we bound the integrand using the upper bound given in Theorem \ref{MainTheorem} with $x=\varepsilon\sqrt n$. Then, 
\begin{equation} \label{ProofMainCor3}
	\E[d_{\textsf{H}}(\hat G^*,G_\mu)\mathds 1_{\hat G\neq\emptyset}]=O\left(n^{-1/2}\right),
\end{equation}
with multiplicative constants that depend on $d,\varepsilon,r,R$ and $L$ only. Together with \eqref{ProofMainCor1}, \eqref{ProofMainCor3} yields the desired result. \hfill \textsquare

\paragraph{Proof of Corollary \ref{CorollaryStoch}}

It is enough to prove that if $\mu$ satisfies either Assumption \ref{Ass11} or \ref{Ass12}, then it satisfies both Assumptions \ref{AssA} and \ref{AssB}, for some values of $\varepsilon, L, r$ and $R$. Hence, Theorem \ref{MainTheorem} will apply and yield the desired result. 

Let $X$ be a random variable in $\R^d$ with probability measure $\mu$. If $u\in\mathcal S^{d-1}$, denote by $f_u$ the density of $\langle u,X\rangle$ and by $F_u$ its cumulative distribution function. For $u\in\mathcal S^{d-1}$ and $t\in\R$, let $\DS\phi(u,t)=f_u(t)=\int_{u^\perp}f(tu+v)\diff v$, where the integral is evaluated with respect to the $(d-1)$-dimensional Lebesgue measure on $u^\perp$.

Let $K$ be the support of $\mu$ and let $\DS A=\{(u,t)\in\mathcal S^{d-1}\times\R: (tu+u^\perp)\cap K\neq\emptyset\}$. Note that $\overset{\circ}{A}$ is included in the support of $\phi$. Thus, by Lemma \ref{Contphi}, $\phi$ is continuous on $\overset{\circ}{A}$. 

From now on, we assume that $\mu$ satisfies either Assumption \ref{Ass11} or \ref{Ass12}. For $u\in\mathcal S^{d-1}$, since $F_u$ is continuous, $q_u^\sharp=q_u^\flat$: Denote this value by $q_u$. Let $\DS \alpha_\textsf{max}=\max_{x\in\R^d}D_\mu(x)$.
This quantity is well defined, since $D_\mu$ is upper semicontinuous and quasi-concave (see \citep{Masse2004}). Let $T\in\R^d$ satisfy $D_\mu(T)=\alpha_\textsf{max}$. Since $\mu$ has a connected support and is absolutely continuous with respect to the Lebesgue measure, such a point exists and is unique (see \citep{Masse2004} or Prop. 3.5 in \citep{MasseTheodorescu1994}). Let $\alpha_1$ and $\alpha_2$ be positive numbers such that $\alpha_1<\alpha<\alpha_2<\alpha_\textsf{max}$. For $u\in\mathcal S^{d-1}$, denote by $q_u^{(1)}$ the $(1-\alpha_1)$-quantile of $F_u$ and by $q_u^{(2)}$ the $(1-\alpha_2)$-quantile of $F_u$.
By Lemma \ref{ContQuant28}, $q_u, q_u^{(1)}$ and $q_u^{(2)}$ are continuous functions of $u$. In addition, for all $u\in\mathcal S^{d-1}$, $\DS \langle u,T\rangle < q_u^{(2)} < q_u < q_u^{(1)}$, by definition of the quantiles and by the first part of Lemma \ref{ContQuant28}. Hence, since $\mathcal S^{d-1}$ is compact, there exist positive numbers $r,R$ and $\varepsilon$ with $\varepsilon<r<R$ and such that for all $u\in \mathcal S^{d-1}$,
\begin{equation} \label{inclusions} 
	\langle u,T\rangle + r \leq  q_u^{(1)} \leq q_u-\varepsilon\leq q_u+\varepsilon \leq q_u^{(2)}\leq \langle u,T\rangle + R.
\end{equation}
In particular, the first and last inclusions of \eqref{inclusions} imply that $B(T,r)\subseteq G_\mu\subseteq B(T,R)$. Hence, $\mu$ satisfies Assumption \ref{AssB}. In addition, by a similar argument as in the proof of Lemma \ref{ContQuant28}, the intermediate inclusions show that the compact set $B=\{(u,t):u\in \mathcal S^{d-1}, q_u-\varepsilon\leq t\leq q_u+\varepsilon\}$ is included in the interior of $A$. Hence, by Lemma \ref{Contphi}, $\phi$ is continuous on and positive on $B$, thus, it is bounded from below by a positive constant $L$ on $B$, yielding $F_u(t')-F_u(t)\geq L(t'-t)$, for all $u\in\mathcal S^{d-1}$ and $t,t'\in\R$ such that $q_u-\varepsilon\leq t\leq t'\leq q_u+\varepsilon$. This, together with continuity of $F_u$, for all $u\in\mathcal S^{d-1}$, shows that $\mu$ satisfies Assumption \ref{AssA}, which finally ends the proof of Corollary \ref{CorollaryStoch}. \hfill \textsquare

\paragraph{Proof of Theorem \ref{TheoremComput}:} 

Let $M\geq 1$, $\DS \frac{10\sqrt{5(d+1)}}{L\sqrt n}\leq z <\varepsilon$ and $\DS \delta=1/\sqrt n$. For simplicity, we denote by $q_j=q_{U_j}$ and $\hat q_j=\hat q_{U_j}$, for $j=1,\ldots,M$. Define the events $\DS \mathcal A=\{|\hat q_j-q_j|\leq z, \forall j=1,\ldots,M\}$ and $\DS \mathcal C=\left\{\{U_1,\ldots,U_M\} \mbox{ is a }\delta\mbox{-net of }\mathcal S^{d-1}\right\}$. Let both $\mathcal A$ and $\mathcal C$ hold. Then, by Lemma \ref{Lemma_q_to_h_Disc}, $\DS d_{\textsf H}(\tilde G_M,G_\mu)\leq Cz+4R\delta$, where $\DS C=\frac{R}{r}\frac{1+\varepsilon/r}{1-\varepsilon/r}$. Therefore, by Lemmas \ref{LemmaUnifDevQuantiles} and \ref{LemmaBarB},
\begin{align*} 
	& \PP\left[d_{\textsf{H}}(\tilde G_M,G_\mu) > Cz+\frac{4R}{\sqrt n}\right] \\
	& \hspace{10mm} \leq Ae^{-L^2x^2/2+10\sqrt{5(d+1)}Lx}+6^d\exp\left(-\frac{M}{2d8^{(d-1)/2}n^{d-1}}+(d/2)\log n\right),
\end{align*}
for any $x\in\R$ satisfying $\DS \frac{10\sqrt{5(d+1)}}{L}\leq x<\varepsilon\sqrt n$. \hfill \textsquare

\bibliographystyle{alphaabbr}
\bibliography{Biblio}

\end{document}